\documentclass[reqno,10pt]{amsart}
\usepackage{amssymb,amsmath,amsthm,enumerate,amscd,graphicx,color}
\usepackage{mathtools}% Loads amsmath
\usepackage[usenames,dvipsnames,svgnames,table]{xcolor}
\usepackage{enumitem}
\usepackage[makeroom]{cancel}
\usepackage{pdfsync}
\usepackage{float}
 \usepackage[colorlinks=true, pdfstartview=FitV, linkcolor=blue, citecolor=blue, urlcolor=blue]{hyperref}
\theoremstyle{rem}
\numberwithin{equation}{section}

\newtheorem{theorem}{Theorem}[section]
\newtheorem{corollary}{Corollary}[section]
\newtheorem{lemma}{Lemma}[section]%
\newtheorem{proposition}{Proposition}[section]
\theoremstyle{definition}
\theoremstyle{remark}
\newtheorem{definition}{Definition}[section]
\newtheorem{remark}{Remark}[section]
\newtheorem{example}{Example}[section]

\newcommand{\B}{\mathcal{B} }
\newcommand{\bC}{{\mathbb C}}
\newcommand{\bZ}{{\mathbb Z}}
\newcommand{\bN}{{\mathbb N}}

\newcommand{\fg}{\mathfrak g }

\newcommand{\W}{\mathcal W}

\newcommand*{\st}{\raisebox{-0.05ex}{*}\llap{\raisebox{-1.2ex}{*}}}
\begin{document}

\author{Caleb Fernelius, Natasha Rozhkovskaya}
\address{Department of Mathematics, Kansas State University, Manhattan, KS 66502, USA}
%\email{rozhkovs@math.ksu.edu}

%\newpage
%\msc{,,,,}    %<-------------------
\thanks{}

\begin{abstract}
We describe the action of the infinite-dimensional Lie algebra $\W_{1+\infty}$ and its B-type analogue on Schur and Schur Q-functions, respectively, using  formal distributions framework. We observe  an interesting self-duality property  possessed  by these compact formulas.
\end{abstract}
\title {Generating functions of $\W_{1+\infty}$ action on symmetric functions}  
\maketitle

\section{Introduction}\label{sec1}
The unique non-trivial central extension of the Lie algebra of differential operators on the circle, commonly known as the  $\W_{1+\infty}$ algebra, gained prominence in the 1990's for  its role in   two-dimensional quantum field theory and integrable systems. Several classes of representations of this infinite-dimensional Lie algebra have been well-studied  (e.g. \cite{AFOMO1}, \cite{FKRW}, \cite{KP}, \cite{KR1},  \cite{KR2}, \cite{PSR}  and in particular, review \cite{AFOMO2}).
 Recent interest in the structure of $\W_{1+\infty}$  comes from  the connection of enumerative geometry  to soliton integrable hierarchies. Studies in these areas   revealed that certain generating functions for  intersection numbers on moduli spaces of stable curves provide examples of tau-functions of famous  integrable hierarchies. Classical results state that 
 Schur symmetric functions serve as tau-functions of the KP hierarchy \cite {Sato}, and  Schur  $Q$-functions serve as  tau-functions of the BKP hierarchy \cite{You1},  and other tau-functions of the corresponding hierarchies  can be expressed  as  finite or  infinite linear combinations of  these classical solutions  with  coefficients satisfying certain algebraic constraints. Hence, the  generating functions of intersection numbers that are  tau-functions of the KP or the BKP hierarchies, can be  also  expressed through  these celebrated families \cite{A1,A2,MM,LY1,LY2,LY3}. In this context,    the symmetries of  the KP, KdV, BKP  hierarchies  are often described  through  the action of  $\W_{1+\infty}$  operators and their  analogues.  In  particular, in \cite{LY3} the explicit formulas for  the action of  certain operators from Lie algebra $\W_{1+\infty}$ and its  $B$ type analogue  on the   Schur functions and Schur $Q$- functions  respectively were found, generalizing  particular cases studied in earlier papers, see cf. \cite{LY3} for citations.

In this paper we continue the discussion of the structure  of   $\W_{1+\infty}$ and  its $B$ type analogue in  the language of  formal distributions  that serve as generating functions for  actions of these algebras  on  families of symmetric functions. Due to close connection  of representations of   infinite-dimensional Lie algebras to the areas studied with the tools of vertex operators, this language is natural and provides additional insight in the structure of such actions. We discovered that in terms of generating functions  formulas of action of   $\W_{1+\infty}$ and its $B$ type analogue on Schur symmetric functions  and   Schur $Q$-functions  respectively are surprisingly compact, involve no derivations, but only multiplication operators, and exhibit an interesting duality, as explained in Remarks \ref{rem6_1}, \ref{rem10_1}.  We also show that  $\W_{1+\infty}$ has the structure of a conformal Lie algebra, and provide formulas for the action of generators of   infinite-dimensional Lie algebras $\hat a_\infty$,  $\hat o_\infty$ in terms of formal distributions of operators acting  on symmetric functions.

The  paper is organized  in the following way. 
In Sections \ref{Sec2},   \ref{Sec3} we   review the algebraic structure of the infinite-dimensional Lie algebras  $\hat a_\infty$ and  $\W_{1+\infty}$ in  the language of formal distributions and their realizations through the Clifford algebra of charged free fermions. In Section \ref{Sec2_6} we describe the conformal Lie algebra structure of  $\W_{1+\infty}$. In Section \ref{Sec4}  we collect necessary  facts about symmetric functions.
Sections \ref{Sec5}, \ref{Sec6}  describe the action of algebraic structures on the space of symmetric functions, with  Theorem \ref{thm6_1} as the main result. Particular  examples and some connections with   formulas deduced by other authors are outlined in the same section. 
 In Section \ref{Sec7} we introduce $\W^B_{1+\infty}$ as a subalgebra of  anti-involution  invariants in $\W_{1+\infty}$, and in Section \ref{Sec8}  we express formal distributions  of generators of this subalgebra and the infinite-dimensional Lie algebra  $\hat o_\infty$ through the Clifford algebra of neutral fermions. In Section \ref{Sec9} we review necessary facts about  Schur $Q$-functions, and in Section \ref{Sec10} we prove  Theorem \ref{thm10_1} that provides actions of $\W^B_{1+\infty}$ and $\hat o_\infty$  on families of symmetric functions.

\section{Lie subalgebra $\W_{1+\infty}$ of $  \hat a_\infty$} \label{Sec2}
\subsection{Lie algebra $\hat a_\infty$ }\label{Sec2_1}
For our goals it is convenient to introduce $\W_{1+\infty}$ as a  central extension of a matrix subalgebra, cf. \cite{KR1}. 
Consider  the  Lie algebra   $a_\infty$  of infinite matrices with a  finite number of nonzero diagonals
 \[ a_\infty=\{(a_{ij})_{ i,j\in \bZ} | \quad a_{ij}=0\text{\quad for \quad} |i-j|>>0\}.\] 
Any  element of $ a_\infty$  is a  finite linear combination of matrices of the form
$ \sum_{i\in\bZ}\lambda_iE_{i,i+k}, $
and  the Lie algebra structure 
\[
\left[ \sum \lambda_iE_{i,i+k},\sum \mu_iE_{j,j+k} \right]= \sum (\lambda_i\mu_{i+k}- \lambda_{i+l}\mu_i)E_{i,i+k+l}
\]
comes from the commutation relations of  generators, which are standard matrix units $\{E_{ij}\}_{i,j\in \bZ}$:
 \[
 [E_{ij},E_{kl}]=\delta_{j,k}E_{il} - \delta_{l,i} E_{kj}.
 \]
 The central extension  $\hat a_\infty=a_\infty\oplus \bC C$  is the Lie algebra  
with central element $C$ and  generators  $\{\hat E_{ij}\}_{i,j\in \bZ}$  such that 
\begin{align}\label{eq2_1}
[\hat E_{ij}, \hat  E_{kl}]=\delta_{j,k}\hat  E_{il} - \delta_{l,i} \hat  E_{kj} +\gamma(\hat  E_{ij},\hat  E_{kl})\, C,
\end{align}
and the cocycle $\gamma(\hat  E_{ij},\hat  E_{kl})$ has the values 
$$
\begin{cases}
\gamma(\hat E_{ij}, \hat E_{ji})=1,& \text{if}\quad i\le 0,\quad  j\ge 1,\\
\gamma(\hat E_{ij}, \hat E_{ji})=-1,& \text{if}\quad   i\ge 1, \quad j\le 0,\\
\gamma(\hat E_{ij}, \hat E_{kl})=0, &\text{otherwise.}
 \end{cases}
$$

 \subsection{Formal distributions} \label{Sec2_2}
  Commutation relations  of  $\hat a_\infty$ and other algebraic structures  below can be presented  in the form of relations on generating functions. For this, following \cite{Kac-begin}, \cite{Kac-bomb},  we recall the notion of a formal distribution.  
Let  $W$ be a vector space. A  {\it $W$-\,valued formal distribution} is a bilateral series 
in the  indeterminate $u$ with coefficients  in $W$:
\begin{align*}
a(u)=\sum_{n\in \bZ}a_n u^n,\quad a_n\in W.
\end{align*}
A formal distribution in two and more  indeterminates is defined similarly.
The vector space of  all $W$-valued formal distributions in indeterminate $u$  is denoted as $W[[u,u^{-1}]]$, and 
we also  use the notation  $W[u]$ for the space  of polynomials,  $W[[u]]$ for the space of power series, $W[u, u^{-1}]$ for the space of Laurent polynomials, and $W((u))$   for the space of formal Laurent series. 
 Formal distributions can be added and multiplied by Laurent polynomials, but  in general  multiplication of two  formal distributions is not a well-defined operation  and can be performed only in  special cases.

 The {\it formal delta-distribution}  $\delta(u,z)$ is the $\bC$-valued  formal 
distribution 
\begin{align*}
	\delta(u,z)=\sum _{i\in \bZ}\frac{u^i}{z^{i+1}}
			=i_{u/z}\left(\frac{1}{z-u}\right) + i_{z/u}\left(\frac{1}{u-z}\right), 
\end{align*}
where  
\begin{align*}
	i_{u/z}\left(\frac{1}{z-u}\right)=  \sum_{i\ge 0}\frac{u^{i}}{z^{i+1}}
\end{align*}
is  the expansion  of the  rational function $1/(z-u)$    in the  domain  $|u|<|z|$. 

One has
\begin{align}
	\delta(u,v)= \delta(v,u)&= -\delta(-u,-v),\quad 
	 \partial_u\delta(u,v)=-\partial_v\delta(u,v),\label{eq2_21}\\
	 \delta(u,z)a(z)&= \delta(u,z)a(u)
		 \quad \text{for any formal distribution}\, a(u).\label{eq2_22}			
\end{align}

 Consider the linear space  $\hat a_\infty[[u,u^{-1}, w,w^{-1}]]$   of formal distributions in  two indeterminates with coefficients in   $\hat a_\infty$.  
Combining  generators  $\{\hat E_{ij}\}$ of $\hat a_\infty$  into a formal distribution $ T (u,w)= \sum_{i,j}\hat E_{ij}u^{i-1} w^{-j}$,  
 commutation relations (\ref{eq2_1}) can be presented in  the form  
\begin{align}\label{eq2_4}
	[ T (u,w),  T (v,z)]&=
	\delta(v,w) \left(T (u,z) + i_{z/u}\left(\frac{1}{u-z}\right)C\right)\,   \notag\\
	&-\,\delta(u,z) \left(T (v,w) +i_{w/v}\left(\frac{1}{v-w}\right)C\right),
\end{align}
and $ [ T (u,w), C]=0$.
\begin{remark}
We introduce    notation  for the central part of this commutation relation: 
\begin{align}\label {eq2_5}
	\gamma(u,w,v,z)&=
	\delta(v,w)  i_{z/u}\left(\frac{1}{u-z}\right)
-\,\delta(u,z) i_{w/v}\left(\frac{1}{v-w}\right)\\
&= i_{z/u}\left(\frac{1}{u-z}\right) i_{v/w}\left(\frac{1}{w-v}\right)-  i_{u/z}\left(\frac{1}{z-u}\right) i_{w/v}\left(\frac{1}{v-w}\right).\notag
\end{align}
Clearly, from (\ref{eq2_4}),
$
\gamma(u,w,v,z)= -\gamma(v,z, u,w).
$
\end{remark}
 The following lemma will be useful  for further calculations with $\gamma(u,w,v,z)$.
\begin{lemma}
\begin{align}\label{eq2_71}
\partial_u^k\partial_v^l\gamma(u,w,v,z\, )|_{w=u,z=v}
&=\frac{(-1)^{k} k! \,l!}{(k+l+1)!}\partial^{k+l+1}_v \delta (u,v),
\end{align}
and
\begin{align}\label{eq2_72}
\partial_u^k\partial_v^l\gamma(u,w,v,z\, )|_{w=-u,z=-v}
= \frac{ k! l!}{(k+l+1)!}\partial^{k+l+1}_v \delta (u,-v).
\end{align}
In particular, 
\begin{align*}
\gamma(u,u,v,v\, )=\partial_v \delta (u,v),
\quad 
\gamma(u,-u,v,-v\, )
=\partial_v \delta (u,-v).
\end{align*}

\end{lemma}
\begin{proof}
For  the proof of (\ref{eq2_71}) observe that
\begin{align*}
	&\partial_u^k\partial_v^l\gamma(u,w,v,z))|_{w=u,z=v}\\
	&=\left.
		\left( i_{z/u}\left(\frac{ (-1)^k k! }{(u-z)^{k+1}}\right) i_{v/w}\left(\frac{l!}{(w-v)^{l+1}}\right)
 		-\,  i_{u/z}\left(\frac{ k!}{(z-u)^{k+1}}\right) i_{w/v}\left(\frac{(-1)^l l!}{(v-w)^{l+1}}\right)\right)\right\vert_{{w=u,z=v}}\\
	&= i_{v/u}\left(\frac{ (-1)^k k! l!} {(u-v)^{k+l+2}}\right)-   i_{u/v}\left(\frac{(-1)^l k! l!\,}{(v-u)^{k+l+2}}\right) 
	= \frac{(-1)^{k} k! l!}{(k+l+1)!}\partial^{k+l+1}_v \delta (u,v).
\end{align*}
The proof of (\ref{eq2_72}) follows the same lines. 
\end{proof}

\subsection{Infinite-dimensional Lie algebra  $\W_{1+\infty}$ } \label{Sec2_3}
Introduce  a family of formal distributions 
 $T^{(k)}(u)\in\hat a_\infty[[u, u^{-1}]]$ by 
\begin{align*}
T^{(k)}(u)= \left.\partial^k_u T (u,v)\right\vert_{v=u}.
\end{align*}
Coefficients of the   expansion 
 $
T^{(k)}(u)=
\sum_{r} T^{(k)}_ru^{r-k-1}
$   are given by  
\begin{align}\label{eq2_10}
T^{(k)}_r=\sum_{j}(r+j-1)_k\hat E_{r+j,j},
\end{align}
where $(a)_k= a(a-1)\dots (a-k+1)$
is the falling factorial. 
Using that  $T^{(k)}(u)$ is a coefficient of  the Taylor series expansion 
\begin{align}\label{eq2_11}
T (u+t,u)= e^{t\partial_u} T(u,v)\,  |_{v=u}= \sum_{k=0}^\infty \frac{\partial^k_u T(u,v)\,  |_{v=u} t^k}{k!} = \sum_{k=0}^\infty \frac{T^{(k)}(u)}{k!} t^k,
\end{align}
we deduce the commutation relations  of formal distributions $\{T^{(k)}(u)\}$.
 \begin{proposition}\label{prop2_1}
\begin{align}\label{eq2_12}
	[ T^{(r)} (u),  T^{(k)}(v)]=&
	\sum_{m=0}^{k}  {k\choose m}{\partial_v^{m}\delta(v,u)\,T^{(r+k-m)}(v)}\notag\\
	&-\sum_{m=0}^{r} {r\choose m}{\partial_u^{m}} \delta(u,v)\,T^{(r+k-m)}(u) \notag\\
	&+\frac{(-1)^{r} r! k!}{(k+r+1)!}\partial^{k+r+1}_v\delta (u,v) C.
\end{align}
\end{proposition}
 \begin{proof}
 Using  (\ref{eq2_22}), (\ref{eq2_11})  we write Taylor expansions of the terms of the commutator (\ref{eq2_4}):
\begin{align}
 	\delta(v+s,u) T (u+t,v) &= \delta(v+s,u) T (v+t+s,v) = \delta(v+s,u) \sum_{k=0}^\infty \frac{T^{(k)}(v)}{k!} (s+t)^k \notag \\
	&= \sum_{r=0}^\infty \sum_{k=0}^\infty   \frac{s^k t^r}{k! r!} \sum_{m=0}^{k} {k\choose m}{\partial_v^{m}\delta(v,u)\,T^{(r+k-m)}(v)}.\label{eq2_12b}
\end{align}
Similarly, with 
\begin{align}
	\delta(u+t,v) T (v+s,u) 
	&= \sum_{r=0}^\infty \sum_{k=0}^\infty   \frac{s^k t^r}{k! r!} \sum_{m=0}^{r} {r\choose m}{\partial_u^{m}\delta(u,v)\,T^{(r+k-m)}(u)},\label{eq2_13}
\end{align}
and applying  (\ref{eq2_71}) in the Taylor expansion 
\begin{align}\label{eq2_14}
\gamma(u+t,u,v+t,v)
	&=\sum_{k=0}^{\infty}\sum_{r=0}^{\infty}    \frac{s^k t^r}{k! r!}  \partial_u^k  \partial_v^k  \gamma(u,w,v,z) \, |_{w=u, z=v},
\end{align}
  we  substitute (\ref{eq2_12b}), (\ref{eq2_13}), (\ref{eq2_14})  into    (\ref{eq2_4})  to get   commutation relations $
[ T (u+t,u),  T (v+s,v)]$.  Proposition \ref{prop2_1} follows from (\ref{eq2_11}).
\end{proof}
Note that Proposition \ref{prop2_1} implies that the linear span of 
   coefficients  of $\{T^{(k)}_r\}_{r\in \bZ, k\in \bZ_{\ge 0}}$ and the central element $C$  form   a Lie subalgebra of 
$\hat  a_\infty$, commonly denoted as  $\W_{1+\infty}$. 
   Commutation relations in $\W_{1+\infty}$ are described by the commutation relations (\ref{eq2_12}) of  generating functions $\{T^{(k)}(u) \}_{k\in \bZ_{\ge 0}}$.

\subsection{$\W_{1+\infty}$ as the  central extension of the  algebra of differential operators}\label{Sec2_4}
The Lie algebra  $\W_{1+\infty}$  is isomorphic to the unique non-trivial central extension $\hat D$ of the  Lie algebra  $D$ of complex regular differential operators on a circle $\bC^{\times}$, see e.g.  \cite {KR1},  
 \cite {VDL} for more details. Lie algebra $D$ has a natural basis $\{t^r(t\partial_t)^k|r\in \bZ, k\in \bZ_{\ge 0}\}$, where $t$ is the variable on $\bC^{\times}$. The natural action of $D$ on the space $t^s\bC[t,t^{-1}]$, where $s\in \bZ$,  defines a homomorphism of Lie algebras $\varphi_s: D\to a_\infty$, which can be  extended to the homomorphism of central extensions $\varphi_s: \hat D\to \hat a_\infty$. 
 On the basis elements the homomorphism  is defined as
\[
\varphi_s(t^r (t\partial_t)^k)= \sum _{j\in \bZ} (-j+s)^k \hat E_{j-r,j}.
\]

Then from (\ref{eq2_10}), for $l\in \bZ_{\ge 0}$, $k\in \bZ$,
$
T^{(k)}_{r}= \varphi_s(t^{-r}\,(-t\partial_t +s-r-1)_{k})$.
Since the highest degree term  of  $t^{-r}\,(-t\partial_t +s-r-1)_{k}$ is $(-1)^k t^{-r+k}\partial_t^k$, 
the set $\{T^{(k)}_{r}\}_{ k\in \bZ_{\ge 0}, r\in \bZ.}$ forms a linear basis of 
 $\W_{1+\infty}\simeq \hat D$.

\subsection{Definition of conformal Lie algebra}\label{Sec2_5}
 Proposition \ref{prop2_1} implies that $\W_{1+\infty}$  defines a  conformal Lie algebra. We review   the  definition of conformal Lie algebra and the $\lambda$-bracket following  \cite{Kac-bomb},\cite{BKV}.

\begin{definition} 
Let   $A$ be a $\bC[\partial]$-module  with the $\lambda$-bracket operation $[a_\lambda b]$,  which is a linear map $ A\otimes A\to A[\lambda]= \bC[\lambda]\otimes A$
 satisfying  the axioms
\begin{enumerate}
\item   (Skewcommutativity) $[ a_\lambda b]= -[b_{-\lambda-\partial} a]$;
\item(Sesquilinearity)
$[\partial a_\lambda b]= -\lambda [a_\lambda b]$,
$[a_\lambda \partial b]= (\partial +\lambda) [a_\lambda b]$ 
\item   (Jacobi identity) $[ a_\lambda [b_\mu c]]=[ [a_\lambda b]_{\lambda+\mu} c]+[ b_\mu[a_\lambda c]]$.
\end{enumerate}
Then $A$ is called a {\it  conformal   Lie algebra}. 
\end{definition}
Let $\fg$ be a Lie algebra.
 Recall  \cite{Kac-begin}, \cite{Kac-bomb} that 
a pair of  $\fg$-valued  formal distributions $a(u)= \sum _{n\in \bZ} a_n u^{-n-1}$ and $b(v)= \sum _{n\in \bZ} b_n v^{-n-1}$  is local  if
\[
(u-v)^N[a(u),b(v)]=0 \quad \text{for} \quad N>>0.
\]
Locality property is equivalent to the existence of  the Operator Product Expansion  (OPE), which means that the commutator of two formal distributions can be expressed through derivatives of formal delta distribution as
\begin{align}\label{OPE}
[a (u),  b(v)]=&\sum_{j=0}^{N-1} c_j\, \partial^j_{v}\delta(u,v)/j!,
\end{align}
where  the coefficient $c_j$ is denoted as $(a(v)_{(j)}b(v))$ and  is called the  $j$-th product of $a$ and $b$. One has  
\begin{align*}
(a(v)_{(j)}b(v)) = \text{Res}_u (u-v)^j[a(u), b(v)].
\end{align*}

Many examples of conformal Lie algebras arise  from families of pairwise  local  distributions, where  the  OPE (\ref{OPE})    defines their  lambda-bracket 
\begin{align*}
[a (u),  b(v)]=&\sum_{j=0}^{N-1} (a(v)_{(j)}b(v))\frac{ \lambda^j}{j!}.
\end{align*}
\subsection{Conformal Lie algebra structure of  $\W_{1+\infty}$ }\label{Sec2_6}
 An adjusted version of commutation relations  (\ref{eq2_12})  implies that every pair $(T^{(k)}(u),\,T^{(l)}(v) )$ is local, and that  the family of formal distributions $\{T^{(k)}(u) \}_{k\in \bZ_{\ge 0}}$ and their derivatives
  generates  a conformal Lie algebra.
 \begin{proposition}\label{prop_2}
 Let $\mathcal T$ be the minimal  $\partial_u$-invariant  subspace of $\hat a_\infty[[u, u^{-1}]]$ that  contains   all  $\{T^{(k)}(u)\}$. Then  $\mathcal T$ has the structure of a conformal  Lie algebra  with the 
 $\lambda$-bracket
\begin{align}\label{eq2_18}
[ T^{(r)} \, _\lambda\,   T^{(k)}]=&
\quad \sum_{m=0}^{k}  {k\choose m}  \lambda^m   {T^{(r+k-m)}}- 
\sum_{m=0}^{r} (-1)^{m}{r\choose {m}}(\lambda+ \partial)^{m}T^{(r+k-m)}  \\
&+\frac{(-1)^{r} r! k!}{(k+r+1)!}\lambda^{k+r+1}\, C.\notag 
\end{align}
\end{proposition}
\begin{proof}

Note that  (\ref{eq2_12}) has a nice symmetry in variables $u$, $v$, but it  is not in the  OPE form (\ref{OPE})  that  would  immediately  define  a conformal Lie algebra  structure, as  it   involves not only derivatives  $\partial_v^{m}\delta(u,v)$, but also    $\partial_u^{m}\delta(u,v)$.  We convert  (\ref{eq2_12})  to the OPE form using the following lemma.

\begin{lemma}\label{lem2_2}
For any formal distribution $a(u)$,
	\begin{align*}
		\partial_u^m \delta(u,v)\,  a(u)
		&=\sum_{p=0}^m (-1)^m{m\choose p}\partial^{m-p}_v\delta(u,v)\,  \partial _v^p a(v).
	\end{align*}
\end{lemma}
\begin{proof}
is by induction on $k$. The base of the induction  for $k=0$ is given by (\ref{eq2_22}).
Assume that for all  values less or equal to  $k$  and for any formal  distribution $a(u)$ the  statement holds.
Note that  from  (\ref{eq2_21}),
\begin{align}\label{eq2_19}
\partial_u^{l} \partial_v^{m}\delta(v,u)= \partial_v^{m}(\partial_u^{l}  \delta(u,v))= (-1)^l
\partial_v^{m+l}  \delta(u,v).
\end{align}
 By  the induction hypothesis,  the Leibniz rule and (\ref{eq2_19}),
\begin{align*}
&\partial_u^{k+1} \delta(u,v) \, a(u)=\partial_u \left(\partial_u^{k} \delta(u,v) \, a(u)\right) -\partial_u^{k} \delta(u,v)\,  \partial _u a(u) 
\\
&=\partial_u \left(
\sum_{p=0}^k(-1)^k{k\choose p}\partial^{k-p}_v\delta(u,v)\, \partial _v^p a(v) \right) -
\sum_{p=0}^k(-1)^{k}{k\choose p}\partial^{k-p}_v\delta(u,v)\,\partial _v^{p+1} a(v)\\
&=
\sum_{p=0}^k (-1)^{k+1}{k\choose p}\partial^{k-p+1}_v\delta(u,v) \partial _v^p a(v)\ +
\sum_{p=1}^{k+1} (-1)^{k+1}{k\choose p-1}\partial^{k-p+1}_v\delta(u,v)\, \partial _v^{p} a(v)\\
&=
\sum_{p=1}^k (-1)^{k+1}\left({k\choose p}+ {k\choose p-1}\right)\partial^{k-p+1}_v\delta(u,v)\partial _v^p a(v) + (-1)^{k+1}\partial^{k+1}_v\delta(u,v)a(u) \, +\,\\
&
+ (-1)^{k+1} \delta(u,v)\partial _v^{k+1} a(v)
=
\sum_{p=0}^{k+1} (-1)^{k+1}{k+1\choose p}\partial^{k+1-p}_v\delta(u,v) \partial _v^p a(u).
\end{align*}
\end{proof}
By  Lemma \ref{lem2_2},  the term   ${\partial_u^{m}} \delta(u,v)\,T^{(r+k-m)}(u)$  in  (\ref{eq2_12})
contributes  to $\lambda$-bracket as 
\[\sum_{p=0}^m (-1)^m{m\choose p}\lambda^{m-p}  \partial _v^p T^{(r+k-m)}(v)= (\lambda+\partial)^{m}  T^{(r+k-m)}(v), \]
which proves  relation (\ref{eq2_18}).
Commutation relation (\ref{eq2_18})  implies that formal distributions  $\{T^{(k)}(u)\}$ are pairwise local, and that the $j$-th products  $T^{(k)}\,_{(j)}T^{(l)}$ are linear combinations of derivatives  of  $\{T^{(k)}(u)\}$.  Hence, $\mathcal T $  is closed under all $j$-th products, and it is a conformal Lie  algebra.

\end{proof}

\subsection{Heisenberg  and Virasoro algebras}  \label{Sec2_7}
 One of the advantageous features of the $\lambda$-bracket form  (\ref{eq2_18}) of commutation relations is that it  effortlessly  reveals the  presence of  several  important  algebraic structures. 
\begin{example}
For $r=k$  relation  (\ref{eq2_18})   takes the  form 
\begin{align}
[ T^{(k)} \, _\lambda\,   T^{(k)}]=
\quad \sum_{m=0}^{k}  {k\choose m}  (\lambda^m  -(-1)^{m}(\lambda+ \partial)^{m})   {T^{(2k-m)}}
+\frac{(-1)^{k} k! k!}{(2k+1)!}\lambda^{2k+1}\, C.\notag 
\end{align}
For   $r=k=0$  the  formula further  reduces the defining   relations of the  Heisenberg algebra: 
\begin{align*}
[ T^{(0)} \,_ \lambda  T^{(0)}]=
\lambda C, 
\end{align*}
since it encodes the commutator of formal distributions  $ [ T^{(0)}(u) , T^{(0)}(v)]=\partial_v\delta(u,v) C$.

For $r=k=1$  we obtain the $\lambda$-bracket  
\begin{align}\label{eq2_20}
[ T^{(1)} \,_\lambda \, T^{(1)}]=&
\partial T^{(1)}+ 2\lambda T^{(1)} -\frac{1}{6} \lambda^3C,
\end{align}
which is the  OPE of  a   Virasoro formal distribution: 
\begin{align*}
[ T^{(1)} (u), T^{(1)}(v)]=&
\partial T^{(1)}(v)\delta(u,v)+ 2 \partial_v\delta(u,v)T^{(1)}(v) -\frac{1}{6} \partial^3_v\delta(u,v)C.
\end{align*}
\begin{remark}
Note  that original relation  (\ref{eq2_12})  is a  {\it symmetric} in $u$ and $v$ version of the Virasoro commutation relation:
\begin{align*}
[T^{(1)}_m, T^{(1)}_n]= \partial_v\delta(u,v)T^{(1)}(v) -  \partial_v\delta(u,v)T^{(1)}(u) -\frac{1}{6} \partial^3_v\delta(u,v)C.
\end{align*}
\end{remark}

Similarly, for $r=k=2$  one gets
\begin{align*}
[ T^{(2)} \,_\lambda \, T^{(2)}]=&
\partial ( 2 T^{(3)}- \partial T^{(2)}) +\lambda( 2 T^{(3)}- \partial T^{(2)})  + \frac{1}{30}\lambda^5C, 
\end{align*}
 and for $k=r=3$  the  $\lambda$-bracket of the  OPE looks like 
\begin{align*}
[ T^{(3)} \,_\lambda \, T^{(3)}]=&
3 \partial  T^{(5)}-3 \partial ^2 T^{(4)}+ \partial ^3 T^{(3)}
+\,  \lambda\, (6 T^{(5)}-6 \partial  T^{(4)}+3 \partial ^2 T^{(3)})\\
&+\,3\lambda^2\, \partial^2 T^{(3)}
+ 2 \lambda^3 T^{(3)}
 - \frac{1}{140}\lambda^7C.
\end{align*}
\end{example}

\begin{example}
 Substitution  $r=1$   in (\ref{eq2_18}) and the sesquilinearity  property describes  the action of the Virasoro algebra on the generators of the algebra $\W_{1+\infty}$:
  \begin{align*}
[ T^{(1)}\, _\lambda  T^{(k)}]=\partial T^{(k)} + (k+1) \lambda T^{(k)} +
\sum_{m=2}^{k}  \lambda^m  {k\choose m}T^{(k+1-m)} -\frac{1}{(k+1)(k+2)} \lambda^{k+2}C.\\
\end{align*}
In particular, 
 \begin{align*}
	[ T^{(1)}\, _\lambda  T^{(0)}]
	=(\partial+\lambda) T^{(0)} +\frac{ 1}{2} \lambda^2C, 
\end{align*}
 which illustrates together with  (\ref{eq2_20})  that coefficients of  formal distributions $T^{(0)}$, $T^{(1)}$,  and  of their derivatives   generate the   Heisenberg-Virasoro subalgebra in    $\W_{1+\infty}$, and that  $T^{(0)}$ and $T^{(1)}$  are  conformal vectors of the Virasoro algebra  action   on $\W_{1+\infty}$
 with  conformal  weights $1$ and $2$ accordingly.  
\end{example}

\section{Presentation  of $\hat a_\infty$ and  $\W_{1+\infty}$ through fermions}\label{Sec3}
\subsection{Clifford  algebra of charged free fermions}\label{Sec3_1}
It is well-known that the Lie  algebra $\hat a_\infty$ and its subalgebra $\W_{1+\infty}$ can be realized through the Clifford  algebra of charged free fermions. 
We review  this presentation in the language of formal distributions,   \cite{Kac-begin},\cite{ Kac-bomb}.
\begin{definition}
The Clifford algebra of  {\it  charged free fermions} is generated by elements  $\{\psi^{\pm}_{i}\}_{i\in \bZ+1/2}$ 
that satisfy relations 
\begin{align*}
	\psi_k^\pm\psi_l^\pm +\psi_l^\pm\psi_k^\pm=0,\quad 
	\psi_k^+\psi_l^- +\psi_l^-\psi_k^+=\delta_{k, -l}, \quad k,l\in \bZ+1/2.
\end{align*}
Collecting generators  as coefficients of     formal distributions
$\psi ^\pm(u)= \sum_{i\in \bZ+1/2}\psi^{\pm}_i u^{-i-1/2}$
 brings  commutation relations  into the form
\begin{align}
	\psi^\pm(u)\psi^\pm(v)+ \psi^\pm(v)\psi^\pm(u)&=0,\label{eq3_1}
	\\
	\psi^+(u)\psi^-(v)+ \psi^-(v)\psi^+(u)&=\delta(u,v).\label{eq3_2}
\end{align}
\end{definition}
Given an associative superalgebra  $A=A_{0}\oplus A_{1}$, the parity of a homogeneous element $a\in A_\alpha$, $\alpha \in \bZ/2\bZ$ is defined as $p(a)=\alpha$.
Consider   a pair of  homogeneous formal distributions   $a(z)=\sum_{n\in \bZ} a_{(n)}z^{-n-1}$  and  $b(w)=\sum_{n\in \bZ} b_{(n)}w^{-n-1}$   with  coefficients in $A$.
Recall  the definition of the normal ordered product $ :a(z)b(w):$ (see e.g.  \cite{Kac-begin}, \cite{Kac-bomb}).
\begin{definition}\label{def3_1}
The {\it normal ordered product} $:a(z)b(w):$   is a formal distribution
\[
 :a(z)b(w):= a(z)_+b(w)+  (-1)^{p(a) p(b)} b(w)a(z)_-,
\]
where  
$a(z)_+=\sum_{n\le -1} a_{(n)}z^{-n-1}$,   $a(z)_-=\sum_{n\ge0}  a_{(n)}z^{-n-1}$,
and  $p(a)$, $p(b)$ are parities of coefficients of $a(z)$ and $ b(z)$ respectively. 
 \end{definition}
Relations between normal ordered products  and regular products of charged  free  fermions can be established by direct calculations.
\begin{lemma}\label{lem3_1}
	\begin{align*}
 		 & :\psi^\pm(u)\psi^\pm(v):=\psi^\pm(u)   \psi^\pm(v),\\
		 & :\psi^+(u)\psi^-(v): =-  \psi^-(v)\psi^+(u) + i_{u/v} \left(\frac{1}{v-u} \right) =\psi^+(u)   \psi^-(v) -i_{v/u} \left(\frac{1}{u-v} \right).
	\end{align*}
  \end{lemma}
\subsection{Presentation  through charged free fermions}\label{Sec3_2}
\begin{proposition} \label{prop3_1}
The algebra $\hat a_\infty$ can be realized  through  charged free fermions  by  
\begin{align}\label{eq3_3}
T(u,w)\to:\psi^+(u)\psi^-(w):
\quad 
\text{ and}\quad C\to Id.
\end{align}
 Accordingly, the algebra $\W_{1+\infty}$ can be realized as a subalgebra of  the Clifford algebra of charged free fermions  by
 \begin{align}\label{eq3_4}
T^{(k)}(u)\to:\partial^k \psi^+(u)\, \psi^-(u):.
\end{align}
\end{proposition}
\begin{proof}
Using Lemma \ref{lem3_1}, we substitute  the commutator of normal ordered products by commutator of regular  products, rearrange terms to compute their anti-commutators and substitute back by normal ordered products.
Note that all products of formal distributions   in calculations below   are well-defined since  they involve  formal distributions of different variables $u, v, w, z$.
\begin{align*}
[:\psi^+(u)&\psi^-(w):, :\psi^+(v)\psi^-(z):] =
[\psi^+(u)   \psi^-(w), \psi^+(v)   \psi^-(z)]\\
=&\psi^+(u)   \psi^-(w) \psi^+(v)   \psi^-(z) -  \psi^+(v)   \psi^-(z) \psi^+(u)   \psi^-(w)\\
=&- \psi^+(u)   \psi^+(v) \psi^-(w)  \psi^-(z) + \psi^+(u)     \psi^-(z)\delta(w,v) \\
&+  \psi^+(v)  \psi^+(u) \psi^-(z)    \psi^-(w) - \psi^+(v)   \psi^-(w)\delta(z,u) 
\\
=&\quad \psi^+(u)     \psi^-(z)\delta(w,v) - \psi^+(v)   \psi^-(w)\delta(z,u) \\
=&:\psi^+(u)     \psi^-(z):\delta(w,v) - :\psi^+(v)   \psi^-(w):\delta(z,u) + \gamma(u,w,v,z),
\end{align*}
where $\gamma(u,w,v,z)$  is given by (\ref{eq2_5}).
Hence, commutator  $[:\psi^+(u) \psi^-(w):,:\psi^+(v)\psi^-(z):]$  matches the  commutator  (\ref{eq2_4}), and we get  (\ref {eq3_3}).
Since derivation $\partial_u$ is compatible with the normal ordered product operation,   and we obtain   (\ref{eq3_4}).
\end{proof}
\begin{remark} 
Under  identification in Proposition \ref{prop3_1},
\[
\hat E_{ij}\to:\psi^+_{-i+1/2}\psi^-_{j-1/2}:=
\begin{cases}
\quad \psi^+_{-i+1/2}\psi^-_{j-1/2}, & i\ge1,\\
-\psi^-_{j-1/2}\psi^+_{-i+1/2}= \psi^+_{-i+1/2}\psi^-_{j-1/2}- \delta_{i,j}, & i\le 0,
\end{cases} 
\]
 and 
 \begin{align*}
 T^{(k)}_r
 &= \sum_{j\ge 1 -r+} (r+j-1)_{k} \psi^+_{-r-j+1/2}\psi^-_{j-1/2} - \sum_{j\le -r} (r+j-1)_{k}\psi^-_{j-1/2} \psi^+_{-r-j+1/2}.
 \end{align*}
 \end{remark}
\begin{proposition}\label{prop3_2}
\begin{align*}
[T(u,v), \psi ^{+} (z)]&= \psi ^{+}(u) \delta(v,z),\quad \quad \quad\quad 
[T(u,v), \psi ^{-} (z)]= -\psi ^{-}(v) \delta(u,z),
\\
[T^{(k)}(u), \psi ^{+} (z)]&= \partial^k_u\psi ^{+}(u) \, \delta(u,z),\quad \quad\quad 
[T^{(k)} (u), \psi ^{-} (z)]= -\psi ^{-}(u)  \, \partial_u^k \delta(u,z).
\end{align*}
\end{proposition}
\begin{proof}
These commutation relations immediately follow from  Lemma  \ref{lem3_1}  with  substitution  of normal ordered products by  regular  products. For example, 
\begin{align*}
[T(u,v), \psi ^{\pm} (z)]&= [:\psi^+(u)\psi^-(v):,\psi ^{\pm} (z)]=  [\psi^+(u)\psi^-(v),\psi ^{\pm} (z)]\\
&=  \psi^+(u)\psi^-(v)\psi ^{\pm} (z) -  \psi ^{\pm} (z) \psi^+(u)\psi^-(v),
\end{align*}
and then  relations (\ref{eq3_1}, \ref{eq3_2}) are applied. 
\end{proof}

\begin{remark}
The last two   relations of Proposition  (\ref{prop3_2})  describe the action of conformal Lie algebra of $\W_{1+\infty}$ on the space of charged free fermions:
 \begin{align*}
[T^{(k)}_\lambda \psi ^{+} ]= \partial^k\psi ^{+},\quad \quad 
[T^{(k)}_\lambda \psi ^{-} ]= (-1)^{k+1}(\lambda+\partial)^k\psi ^{-}.
\end{align*}
\end{remark}

\section{Review of properties of  symmetric functions}\label{Sec4}
\subsection{On the boson-ferimon correspondence}\label{Sec4_1}
In this section we construct the action of 
Lie algebras  $\hat a_\infty$ and  $\W_{1+\infty}$  on a vector space that contains  countably many copies of the ring of symmetric functions. 
 This construction is based on the   discussed in Section \ref{Sec3} presentation of  Lie algebras  through charged free fermions
and on the {{\it boson-fermion correspondence}}.
 We only  briefly outline   the main idea of this important  isomorphism of modules of infinite-dimensional algebraic structures,  referring  
for more details to many books and papers, such as e.g.\cite{Kac-begin},  \cite {Kac-bomb}, \cite{MJD}. For our purposes we will need only the
 resulting from the boson-fermion correspondence
formulas (\ref{eq5_3}), (\ref{eq5_4}) of action of charged free fermions on the ring of symmetric functions. 

The boson-fermion correspondence establishes an isomorphism of two infinite-dimensional vector spaces as modules over several important algebraic structures. 
The first  of the two vector spaces  is  the space of semi-infinite wedge forms, usually called the fermionic Fock space. The second one is  the space of
  countably many copies of a polynomial ring of infinitely many variables, called the bosonic Fock space. It turns out  that 
   the Clifford algebra  of  charged free fermions, Heisenberg algebra, Virasoro algebra  act  on both spaces in an equivalent way, and  moreover, the  actions of these algebraic structures  are closely related to each other. 

 In particular,   the Clifford algebra of charged free fermions  has a natural action  on  the fermionic Fock space of semi-infinte  wedge forms \cite{Kac-bomb}. The boson-fermion correspondence transports this action into the  action  on  the boson Fock space. An important  step of this  transition  
  is the identification of  each  polynomial ring that constitute graded components of the boson Fock space with the ring of symmetric functions. 
This way  charged free fermions act  on the space   of countably many copies  of the ring of symmetric functions. We recall the explicit formulas of
   this action below  in  Section \ref{Sec5_1}.
   
Presentations of   $\hat a_\infty$ and  $\W_{1+\infty}$  through  charged free fermions provides actions of
 these two Lie algebras on the space of countably many copies  of the ring of symmetric functions. 
 Our goal  is to describe properties of these actions in terms of generating functions, which is done in   Section \ref{Sec6}.
  In the rest of this section we review  the necessary facts about symmetric functions. 

\subsection{Symmetric functions}\label{Sec4_2}
 We review properties  of symmetric functions   \cite{Md}, \cite{Stan} in  the   setup and notations   similar to   \cite{JR-genA}, \cite{KLR}, \cite{NC}, \cite{NRQ}.
Consider the  algebra of formal power series $\bC[[{\bf x}]]=\bC[[x_1,x_2,\dots]]$. 
Let $\lambda=(\lambda_1\ge \dots\ge \lambda_l> 0)$ be a partition  of length $l$.
The {\it monomial symmetric function} 
is  a formal series
\[
m_\lambda=\sum_{(i_1,\dots,i_l)\in \mathbb N^l} x^{\lambda_1}_{i_1}\dots x^{\lambda_l}_{i_l}.
\]
Let  $\Lambda$ be the  subalgebra  of $\bC[[{\bf x}]]$  spanned as a vector space  by all  {monomial symmetric functions}. It is  called the  {\it algebra of symmetric functions}.  Note that elements  of $\Lambda$ are invariant with respect to any permutation of  a finite number of  indeterminates 
 $x_1, x_2, \, \dots$. The following  families of symmetric functions play important role in our study. 

For a partition $\lambda=(\lambda_1\ge \dots\ge \lambda_l> 0)$,  {\it Schur symmetric   function} $s_\lambda$ is defined as 
\begin{align}\label{eq4_1}
s_\lambda(x_1,x_2,\dots ) =\sum_{T} {\bf x}^{T},
\end{align}
where the sum is over all semi-standard tableaux of shape $\lambda$.

Complete symmetric functions  $h_k= s_{(k)}$,  elementary symmetric functions $e_k= s_{(1^k)}$,  and power sums $p_k$ are defined 
 by the formulas
\begin{align}
h_k&=\sum_{i_1\le i_2\le \dots \le i_k} x_{i_1} x_{i_2}\dots x_{i_k},\quad \label{eq4_2}\\
 e_k&=\sum_{i_1< i_2<\dots < i_k} x_{i_1} x_{i_2}\dots x_{i_k},\quad\quad  \label{eq4_3}\\
  p_k&= \sum_{i\ge 1} x_i^k. \label{eq4_4}
\end{align}
Set
$h_{-k}(x_1,x_2\dots)=e_{-k}(x_1,x_2\dots)=p_{-k}(x_1,x_2\dots)=0
$
 for 
$ k\in \bN$  and $h_0=e_0=p_0=1$.

The algebra $\Lambda$  is  a polynomial algebra in each of these three families of generators: 
\[
 \Lambda=\bC[h_1, h_2,\dots]=\bC[e_1, e_2,\dots]=\bC[p_1,p_2,\dots].
\]
Schur symmetric functions $\{s_\lambda\}$ labeled by all partitions  form a  linear basis of $\Lambda$.
They can be also expressed through complete symmetric functions  by the {\it Jacobi\,-\,Trudi identity}
\begin{align}\label{eq4_5}
s_\lambda=\det[h_{\lambda_i-i+j}]_{1\le i,j\le l}.
\end{align}
We will use (\ref{eq4_5}) as the extension of the  definition of $s_\lambda$  assuming that $\lambda$ can be any integer vector.
Then for any  $\lambda\in \bZ^l$ either $s_\lambda$  coincides  up to a sign with a Schur symmetric function (\ref{eq4_1}) labeled by a partition, or equals zero.

Generating series  of the above families of symmetric functions are  elements of  $ \Lambda[[u]]$ that are  given by 
\begin{align}\label{eq4_6}
	H(u)=\sum_{k \ge 0}  {h_k}{u^k}&=\prod_{i\in \bN} \frac{1}{1-x_iu},\quad 
	E(u)=\sum_{k\ge 0}{ e_k}{ u^k}=\prod_{i\in \bN} {(1+x_iu)}.
\end{align}
We also denote 	$P(u)= \sum_{k\ge 1}{ p_k}{ u^{k-1}}$.

 Most of the time we will not recall dependence   (\ref{eq4_2})-(\ref{eq4_4}) of symmetric functions on the original  variables of $x_i$'s, 
 but  treat them as a polynomials in one of   families  of generators,   $\{p_i\}$, or $\{e_i\}$, or  $\{h_i\}$. 
In some cases below  we  still  need to specify   the set  of   indeterminates $(x_1, x_2\dots )$ as   in definitions (\ref{eq4_2}) - (\ref{eq4_4}). 
Then we  write
$h_k=h_k(x_1,x_2,\dots)$,  $ e_k=e_k(x_1,x_2,\dots)$, $p_k=p_k(x_1,x_2,\dots)$,  $H(u; x_1, x_2,\dots)$, $E(u; x_1, x_2,\dots)$,  $P(u; x_1, x_2,\dots)$, etc.

There is a natural  scalar product on   $\Lambda$,  where  the set of Schur symmetric functions $\{s_\lambda\}$ labeled by partitions $\lambda$ form an orthonormal basis,
$
 	\langle s_\lambda,s_\mu \rangle=\delta_{\lambda, \mu}.
$
Then, for any  linear operator  acting on the  vector space $\Lambda$, one can define the corresponding adjoint operator. In particular,  any symmetric function  $f\in \Lambda $  defines an operator of multiplication
$
f:g\mapsto fg$  for any $g\in \Lambda$.
The corresponding   adjoint operator $f^\perp$  is defined by the standard rule
$\langle f^\perp g_1, g_2\rangle =\langle g_1,fg_2 \rangle $ for all   $g_1,g_2\in \Lambda$.

It is known  \cite[I.5  Example 3]{Md}  that
$
p_n^\perp= n\frac{\partial}{\partial p_n}.
$
Since any element  $f\in \Lambda$  can be expressed as a polynomial  function of power sums 
\begin{align*}
f= F(p_1,p_2, p_3, \dots),
\end{align*}
the corresponding adjoint operator $f^\perp$ is a polynomial differential operator with constant coefficients
\begin{align*}
f^\perp= F(\partial/\partial p_1,2\partial/\partial p_2, 3\partial/\partial p_3,\dots).
\end{align*}
In particular, $e_k$ and $h_k$ are  homogeneous polynomials of degree $k$ in  $(p_1,p_2, p_3,\dots )$, so
 the adjoint operators $e^\perp_k$ and $h^\perp_k$ are  homogeneous polynomials of degree $k$ in $(\partial/\partial p_1,2\partial/\partial p_2, \dots)$, which implies
  the following statement. 
\begin{lemma}\label{lem4_1}
For any symmetric function $f\in \Lambda$   there exists  a positive  integer   $N= N(f)$, such that
 \[e^\perp_l(f)=0\quad \text{and}\quad h^\perp_l(f)=0\quad \text{ for all\quad   $l\ge N$}.
 \] 
 \end{lemma}
 Lemma \ref{lem4_1} guarantees that formal distributions of operators that we use later in this text act  locally finitely on the boson Fock space and that the products
of these formal distributions of operators  are well-defined.

Denote by $\mathcal D$  the algebra  of differential operators  acting on $\Lambda=\bC[p_1,p_2,\dots]$, which consists of finite sums  
\[
\sum_{i_1,\dots, i_m} F_{i_1, \dots i_m}(p_1, p_2,\dots) {\partial _{p_1}^{i_1}}\dots  {\partial _{p_m}^{i_m}},
\] 
where coefficients  $F_{i_1, \dots i_m}(p_1, p_2,\dots)$ are  polynomials in $(p_1, p_2,\dots)$.
Then operators of multiplication $p_n, h_n, e_n$, their adjoints $p_n^\perp, h_n^\perp, e_n^\perp$ along with 
 their products  are elements of   $\mathcal D$.

We will use the same notation for  the generating series of the corresponding multiplication operators:
$H(u),E(u), P(u)\in \mathcal D[[u]]$. Similarly, we define $E^\perp(u), H^\perp(u),  P^\perp(u)\in \mathcal D[[u^{-1}]]$ as
\begin{align}
  E^\perp(u)= \sum_{k\ge 0}\frac {e^\perp_k} {u^k},\quad H^\perp(u)= \sum_{k\ge 0} \frac{h^\perp_k} {u^k},\label{eq4_7}\\
  \quad P^\perp(u)=\sum_{k\ge  1} \frac{p^\perp_{k}}{u^{k+1}}= \sum_{k\ge  1}\frac{ k{\partial/}{\partial_{p_{k}}} }{u^{k+1}}.\label{eq4_8}
\end{align}

The properties of operators collected in the next proposition  are either well-known or can be easily deduced from well-known  properties of symmetric functions \cite[I.2] {Md}.
  \begin{proposition}\label{prop4_1}
  In  $\mathcal D[[u]]$ (resp. in  $\mathcal D[[ u^{-1}]]$), 
\begin{align}
&H(u) E(-u)=1,\quad 
H^\perp(u) E^{\perp}(-u)=1,\label{HEP1}\\
&P(-u)={\partial_uE(u)}{H(-u)}, \quad   P(u)={\partial_u H(u)}{E(-u)},\label{HEP2}\\
&\partial_u^k H(u)= H(u)(\partial_u+P(u))^{k-1} (\,P(u)\,),\label{HEP3}\\  
&\partial_u^k E(u)= E(u)(\partial_u+P(-u))^{k-1}( \,P(-u)),\label{HEP4}\\
&{ P^\perp(-u)=- \partial_u E^\perp(u) \,\, H^\perp(-u),\quad }
P^\perp(u)= - \partial_u H^\perp(u)\, \, E^\perp(-u),\label{HEP7}\\
& H(u)= exp\left(\sum_{k\ge1} \frac{p_k}{k}{u^k}\right),\quad 
  E(u)= exp\left(-\sum_{k\ge1} \frac{(-1)^{k} p_k}{k}{u^k}\right), \label{HEP5}
\\
&H^\perp(u)= exp\left( \sum_{k\ge 1}\frac{\partial}{\partial p_k} \frac{1}{u^k}\right),\quad 
E^\perp(u)= exp\left(-\sum_{k\ge 1} {(-1)^k} \frac{\partial}{\partial p_k} \frac {1}{u^k}\right).\label{HEP8}
\end{align}
In (\ref{HEP3}), (\ref{HEP4}),  $(\partial_u+P(\pm u))^{k-1} (\,P(\pm u))$ means that  differential operator $(\partial_u + P(\pm u))^{k-1}$ is applied  to
  formal distribution of operators $P(\pm u)$.
\end{proposition}
  \begin{lemma} \label{lemmaHE}  We have  the following commutation relations in   $\mathcal D[[u^{-1}, v]]$:
  \begin{align}
	\left(
	1-\frac{v}{u}
	\right)E^\perp(u)E(v)&= E(v)E^\perp(u),\label{CR1}\\
\left(1-\frac{v}{u}
\right)H^\perp(u)H(v)&= H(v)H^\perp(u),\label{CR2}\\
	H^\perp(u)E(v)&= \left(
	1+\frac{v}{u}
	\right)E(v)H^\perp(u),\label{CR3}\\
E^\perp(u)H(v)&= \left(
1+\frac{v}{u}
\right)H(v)E^\perp(u),\label{CR4}\\
	P^\perp(u)H(v) &=  H(v)P^\perp(u)-
	\frac{v}{u^2}i_{v/u}\left(1-\frac{v}{u}\right)^{-1}H(v),\label{CR5} 
\\
P^\perp(u) E(v)&=E(v)P^\perp(u)+\frac{v}{u^2} 
i_{v/u} \left( 1+\frac{v}{u} \right)^{-1} E(v).\label{CR6} 
\end{align}
\end{lemma} 
\begin{proof}
Relations (\ref{CR1}) - (\ref{CR4}) can be found 
  \cite[I.5 Example 29] {Md}.  
    For the proof of  (\ref{CR5}), 
  differentiating (\ref{CR2}), and using (\ref{HEP7}) and (\ref{HEP1}),  we get  
  \begin{eqnarray*}
    P^\perp(u)H(v)&=&-E^\perp(-u)\partial_uH^\perp(u)H(v)\\
    &=&-E^\perp(-u)\left(\frac{v}{u^{2}}i_{v/u}\left(1-\frac{v}{u}\right)^{-2}H(v)H^\perp(u)+i_{v/u}\left(1-\frac{v}{u}\right)^{-1}H(v)\partial_uH^\perp(u)\right)\\
    &=&-\frac{v}{u^2}\left(1-\frac{v}{u}\right)i_{v/u}\left(1-\frac{v}{u}\right)^{-2}H(v)E^\perp(-u)H^\perp(u)\\
    & \quad&-\left(1-\frac{v}{u}\right)i_{v/u}\left(1-\frac{v}{u}\right)^{-1}H(v)\partial_uH^\perp(u)E^\perp(-u)\\
    &=&H(v)\left(-\frac{v}{u^2}i_{v/u}\left(1-\frac{v}{u}\right)^{-1}+P^\perp(u)\right).
\end{eqnarray*}
Relation (\ref{CR6}) is proved similarly.
\end{proof}

Finally, we would like to state a simple  corollary of Lemma \ref{lemmaHE} that  will be used in Sections \ref{Sec6} and \ref{Sec10}.
\begin{corollary}\label{cor4_1}
  In   $\mathcal D[[u^{-1}, u_1,\dots, u_l]]$
\begin{align*}
  E^\perp(-u)\prod_{i=1}^{l} H(u_i)&= \prod_{s=1}^l \left(1-\frac{u_i}{u}\right) \prod_{i=1}^{l} H(u_i)E^\perp(-u),\\
 H^\perp(u)  \prod_{i=1}^{l} H(u_i)&=   \prod_{s=1}^l i_{u_i/u}\left(1-\frac{u_i}{u}\right)^{-1} \prod_{i=1}^{l} H(u_i) H^\perp(u),\\
 P^\perp(u)  \prod_{i=1}^{l} H(u_i)&=\prod_{i=1}^{l} H(u_i)  P^\perp(u)-
   \left( \sum_{i=1}^l  \frac{u_i}{u^2}  i_{u_i/u} \left( 1-\frac{u_i}{u} \right)^{-1}\right)\prod_{i=1}^{l} H(u_i).
\end{align*}
\end{corollary}

\section{Action of charged free fermions on the space of  symmetric functions} \label{Sec5}
  
\subsection{Action of charged free fermions  on symmetric functions} \label{Sec5_1}
Consider the ring of symmetric functions $\Lambda$ as the ring of polynomials in power sums   $\Lambda=\bC[p_1,p_2,\dots]$.
The  {\it boson Fock space}  is the graded space  $\B= \bC[ z, z^{-1}]\otimes \Lambda$ that consists of countably many copies of $\Lambda$.
The grading $\B=\oplus_{m\in \bZ} \B^{(m)}$  with  $\B^{(m)}=  z^m\,\Lambda$  is called 
 {\it charge decomposition}.

Let $R(u)$ act  on the  elements of the  form  $z^m f \in \B^{(m)}$,  where $f \in \Lambda$,  $m\in \bZ$, by the rule
\[
R(u) (z^mf)={ z^{m+1}}{u}^{m+1} f.
\]
Then   $R^{-1}(u)$ acts as
\[
R^{-1}(u) (z^mf)={ z}^{m-1} u^{-m}f.
\]
The  meaning of the action  of $R^{\pm 1}(u)$ is the translation  between equivalent  symmetric functions as elements  in  different  graded components of $\B$.

Define  formal distributions $\psi^\pm(u)$  of operators  acting on the  space $\B$  through the action of  $R^{\pm 1}(u)$ and the $\mathcal D$-valued generating series (\ref{eq4_6}), (\ref{eq4_7}):
\begin{align}
\psi^+(u)&=u^{-1}R(u)H(u) E^{\perp}(-u),\label{eq5_1}\\
\psi^-(u) &=R^{-1}(u)E(-u) H^{\perp}(u),\label{eq5_2}
\end{align}
or, in  other words, for any $m\in \bZ$ and any $f\in \Lambda$,
\begin{align}
\psi^+(u)(z^mf)&=z^{m+1}u^{m}H(u) E^{\perp}(-u)(f),\label{eq5_3}\\
\psi^-(u)(z^mf)&=z^{m-1}{u^{-m}}E(-u) H^{\perp}(u)(f). \label{eq5_4}
\end{align}
Let the operators   $\{\psi^{\pm}_{i}\}_{i\in \bZ+1/2}$ be the coefficients of the expansions
 \[
\psi ^\pm(u)= \sum_{i\in \bZ+1/2}\psi^{\pm}_i u^{-i-1/2}.
\]

\begin{proposition}\label{prop5_1} 
Formulas (\ref{eq5_1}),  (\ref{eq5_2}) define quantum fields $\psi^\pm(u)$ of  operators acting on  the space  $\B$ that satisfy relations of charged  free fermions
(\ref{eq3_1}), (\ref{eq3_2}). 
\end{proposition}
For the original proof we refer to \cite{Jing2}. Alternatively, one can deduce the relations
from the definition (\ref{eq5_3}, \ref{eq5_4}) and Lemma \ref{lemmaHE} (see \cite{NC} for this approach).

\begin{remark}\label{rem5_1}
From  (\ref{HEP5}), (\ref{HEP8}) one immediately  gets the bosonic form of   
 $\psi^{\pm}(u)$:
\begin{align*}
 \psi^+(u)&=
 u^{-1}R(u)
 \exp \left(\sum_{n\ge 1}\frac{p_n}{n}{u^n}  \right) \exp \left(-\sum_{n\ge 1}\frac{\partial} {\partial p_n}\frac{1} {u^{n}}\right),
\\
\psi^-(u)&=
R^{-1}(u)
\exp \left(-\sum_{n\ge 1}\frac{p_n}{n} {u^n}  \right) \exp \left(\sum_{n\ge 1}\frac{\partial} {\partial p_n} \frac{1}{u^{n}}\right).
\end{align*}
\end{remark}

 \subsection{Generating function for  Schur  symmetric functions}\label{Sec5_2} 
 Applying commutation relation  (\ref {CR4})   to a product of  multiple operators  of the form  (\ref{eq5_3}),  we get
\begin{align*}
\psi^+(u_l)\dots \psi^+(u_1) \, (z^m f)
=z^{l} \prod_{1\le i<j\le l}\left(1-\frac{u_i}{u_j}\right)\prod_{i=1}^l u_i^{k+i-1}H(u_i) \prod_{i=1}^l E^{\perp}(-u_i) (z^mf).
\end{align*}
In  particular, let 
\[S^{(m)}(u_l,.., u_1)= \psi^+(u_l)\dots \psi^+(u_1) \, (z^m)
\]
be the result of action of  operators $\psi^+(u_l)\dots \psi^+(u_1)$ on the  $m$-th vacuum vector $z^m\in\B^{(m)}$.
Since   $E^{\perp}(-u_i)(z^m)= z^m$, we can write 
\begin{align*}
S^{(m)}(u_l,.., u_1)=
  z^{m+l} u_l^k\dots u_1^k S(u_l,.., u_1),
\end{align*}
where ${S}(u_l,.., u_1)$ is a formal distribution with coefficients in $\Lambda$:
   \begin{align}\label{eq5_6}
{S}(u_l,.., u_1)&=\prod_{1\le i<j\le l}\left( u_j- {u_i}\right)\prod_{i=1}^{l} H(u_i)
\end{align}
In the expansion ${S}(u_l,.., u_1)=\sum_{\alpha\in \bZ} 
S_\alpha u_1^{\alpha_1}\dots u_l^{\alpha_l}$
coefficients  $S_\alpha= \det[ h_{\alpha_i+1-j}] = 0 $ if  $\alpha_i<0 $ for any $i$,  which allows us 
to write this expansion as   the sum  over $ \alpha\in \bZ^{l}_{\ge 0}$ rather than over all $\bZ^l$.
Also $S_\alpha= 0 $ if  $\alpha_i= \alpha_j $ for some  $i\ne j$. 
The transition to standard notations of  Schur symmetric functions
 $s_\lambda$  is given by formulas 
\begin{align*}
 & s_{(\lambda_1,\dots, \lambda_l)} = S_{( \lambda_1+l-1, \lambda_2+l-2,\dots \lambda_l)},\\
 & e_k= s_{(1^k)}=  S_{ (  { \underbrace {l,\,l-1,\,\dots, l-k+1,\,}_{k}  }  \, \,{ \underbrace { l-k-1\,,l-k-2\,,\dots,0\,}_{l-k} } )},\\
  &h_k = s_{(k)}= S_{( k+l-1, \quad {\underbrace{l-k-1,\,l-2,\,\dots,2,\,1,\,0}_{l-1} } )}.
\end{align*}
Every coefficient $S_\alpha$ is  either zero, or coincides  up to a sign with one of   Schur symmetric functions $s_\lambda$, and every Schur symmetric  function  $s_\lambda$ appears as a coefficient in the  expansion of $\mathcal{S}(u_l,.., u_1)$.
 In that sense $\mathcal{S}(u_l,.., u_1)$  can be viewed  as the  generating function for  Schur symmetric functions $\{s_\lambda\}_{\lambda \in \bZ^l}$.

 \begin{remark}\label{rem5_2}
 From (\ref{eq5_6}) and  Cauchy identity (see e.g.\cite{Md} I.4 ,(4.3)) 
in the region $|u_1|<|u_2|<\dots <|u_l|$, 
 \begin{align*}
 {S}(u_l,.., u_1)\prod_{1\le i<j\le l} i_{u_i/u_j}\left( u_j- {u_i}\right)^{-1}
 =\prod_{i=1}^{l} H(u_i)= \sum_\lambda s_\lambda(x_1,x_2,\dots) s_\lambda(u_1, u_2, \dots u_l).
 \end{align*}
 \end{remark}

\subsection{Bosonic  normal ordered products for charged free fermions}\label{Sec5_3}

Following \cite {MJD}, we introduce another version of normal ordered product  $\st\, \st$, defined for polynomial differential   operators acting  on the boson Fock space.
In a product of two formal distributions of such operators this reordering moves  all multiplication   operators  by $p_i$'s  to the left from  all differentiations $\partial/\partial p_i$'s.We we will call it  {\it bosonic  normal ordered product}.
In particular,  for charged free  fermions the {bosonic  normal ordering} is defined as follows:
  \begin{align*}
\st \psi^+(u)\psi^-(v)\st (z^mf) &= z^{m}u^{m}{v^{-m}}H(u) E(-v) E^{\perp}(-u) H^{\perp}(v) (f),\\
\st \psi^+(u)\psi^+(v)\st (z^mf) &=z^{m+2}u^{m+1}{v^{m}}H(u) H(v) E^{\perp}(-u) E^{\perp}(-v) (f),\\
\st \psi^-(u)\psi^-(v)\st (z^mf) &=  z^{m-2}u^{-m+1}{v^{-m}}E(-u) E(-v) H^{\perp}(u) H^{\perp}(v) (f).
\end{align*}
Relations between different products of charged free fermions follow  from  Lemmas \ref{lem3_1} and   \ref{lemmaHE},  and  
  can be summarized as follows:
  \begin{align}
  : \psi^+(u)\psi^-(v):&=i_{u/v} \left(\frac{1}{v-u} \right) -  \psi^-(v)\psi^+(u)  =   i_{u/v} \left( \frac{1}{v-u}\right) ( 1-\st \psi^+(u)\psi^-(v)\st ),\label{t1}\\
    : \psi^\pm(u)\psi^\pm(v):&=   \psi^\pm(u)\psi^\pm(v) =  \left(1-\frac{v}{u} \right) \st \psi^\pm(u)\psi^\pm(v)\st.\label{t2}
\end{align}

We introduce notation  
  $\st T(u,v)\st =\st\, \psi^+(u)\psi^-(v)\,\st$.

  \begin{proposition}\label{prop5_2}
In the bosonic order the generating functions of $\hat a_\infty$ and $\W_{1+\infty}$ have the form
  \begin{align}
  T(u,v)&=  i_{u/v} \left( \frac{1}{v-u}\right) ( 1-\st\, T(u,v) \st), \label{eq5_9}\\
    T^{(k)}(u)&=  \frac{1}{k+1}\partial ^{k+1}_u \st\, T(u,v)\st |_{v=u}.\label{eq5_10}
  \end{align}
    \end{proposition}
    \begin{proof}
Statement (\ref{eq5_9}) immediately follows from (\ref{t1}).
Note that we cannot evaluate $T^{(k)}(u)$  in  the most  straight-forward way as a derivative  of  $-  \psi^-(v)\psi^+(u) + i_{u/v} \left(\frac{1}{v-u} \right)$, as  each of these two  formal distributions has a singularity at $u=v$. However, the bosonic  presentation (\ref{eq5_9}) allows us to overcome this obstacle. Since the linear polynomial   $(v-u)$ is the inverse element for the formal power  series $ i_{u/v} \left( \frac{1}{v-u}\right)$, 
   multiplication of both parts of (\ref{eq5_9}) by $(v-u)$ results in
      \begin{align*}
(v-u)  T(u,v)=  1-\st T(u,v) \st.
  \end{align*}
Application of the operator  $\partial^{k}_u$ to both sides  gives the relation 
          \[
- k\,\partial ^{k-1}_uT(u,v)+ (v-u)\partial ^{k}_u T(u,v)=  -\partial ^{k}_u \,\st T(u,v) \st, 
  \]
  and setting  $v=u$ we get (\ref{eq5_10}).
      \end{proof}
\section{Action of $\hat a_\infty$ and   $\W_{1+\infty}$ on symmetric functions}\label{Sec6}     
\subsection{The first main result}\label{Sec6_1} In this section we describe the action of $\hat a_\infty$ and   $\W_{1+\infty}$ on symmetric functions through the language of formal distributions those coefficients are operators of multiplication and differentiation on the space of symmetric functions. Remark
    \ref{rem6_1}  provides  a comment on  the interesting duality  of  formula (\ref{eq6_4}). 
    To highlight  this duality we  introduce  special  notations for  products of   $\left(1-\frac{u_i}{u}\right)$  and $i_{u_i/u}\left(1-\frac{u_i}{u}\right)^{-1}$.
      Set 
      \[
   {\mathcal E_l} (-u^{-1})=    {E} (-u^{-1};\, u_1, \dots, u_l)= \prod_{s=1}^l \left(1-\frac{u_i}{u}\right),
      \]
     and 
          \[
            {\mathcal H_l(u^{-1})}= \frac{1}{\mathcal E_l(-u^{-1})}  =  {H} (u^{-1};\, u_1, \dots, u_l)=  \prod_{i=1}^{l} i_{u_i/u}\left(1-\frac{u_i}{u}\right)^{-1}.
      \]
     Observe that  
     $ {\mathcal E_l} (-u^{-1})  $ is a polynomial in variables $( u_1,\dots, u_l, u^{-1})$,  and $ {\mathcal H_l(u^{-1})}$
      is a  formal power series  in variables $( u_1,\dots, u_l, u^{-1})$. The introduced notations are  purposefully compatible with (\ref{eq4_6}), providing 
      interpretation of    coefficients of  $ {\mathcal E_l} (-u^{-1})  $ and $ {\mathcal H_l(u^{-1})}$   as   symmetric  polynomials
   $e_k(u_1,\dots, u_l)$ and $h_k(u_1,\dots, u_l)$ respectively.
    
   \begin{theorem}\label{thm6_1} 
     	 \begin{enumerate}
  		  \item 
			Let $f\in \Lambda$. The  action of $\hat a_\infty $  on  $z^mf\in \mathcal B^{(m)}$ is described by 
			\begin{align}\label{eq6_1}%{Tf1}
				T(u,v)(z^m f)=z^{m}   i_{u/v} \left(\frac{1}{v-u}   \right)\left(1-\frac{u^{m}}{v^{m}}E(-v)H(u)E^{\perp}(-u)  H^{\perp}(v) \right)f.
			\end{align}
		\item   For any  $f\in \B^{(0)}$  the action of  $\W_{1+\infty}$ on this element is described by 
   			 \begin{align}\label{eq6_2}%{Tf3}
     				 T^{(k)}(u) \,(f)=  \frac{1}{k+1}\sum_{r+s=k+1}{s\choose{ k+1}} E(-u) \,\partial ^{r}_u H(u)  \,  \,H^{\perp}(u)\, \partial_u^s (E^{\perp }(-u)) \,  (f).
   			 \end{align}

		\item  The action of $\hat a_\infty $ on the basis of Schur symmetric functions is described  by the following relation on generating functions. 
			\begin{align}\label{eq6_3}%{Tf2}
				&T(u,v)\left(S(u_l,.., u_1)^{(m)}\right)\notag\\
				=&\quad  i_{u/v}\left(\frac{1}{v-u}\right) \left(S(u_l,.., u_1)^{(m)} 
				- \frac{1}{zu^{l} v^{m}}{{\mathcal H_l(v^{-1})} } E(-v) 
				S(u, u_l,.., u_1)^{(m)}\right).
			\end{align}
	\item The action  of   $\W_{1+\infty}$ on the basis of Schur symmetric functions is described by the following relation on generating functions. 
 			  \begin{align} \label{eq6_4}%{Tf4}
    				 T^{(k)}(u) \,\left(S(u_l,.., u_1)\right)
	= \frac{\,\partial_u^{k+1}({\mathcal E_l}(-u^{-1}) H(u)\,)}{(k+1) {\mathcal E_l}(-u^{-1}) H(u)} S(u_l,.., u_1).
  			  \end{align}   
		\end{enumerate}
\end{theorem}

\begin{proof}
\begin{enumerate}
\item  Action (\ref{eq6_1}) is immediate consequence of (\ref{eq5_9}).
\item Singularities   of  (\ref{eq6_1}) at  $v=u$ prevent us from  using that formula directly to evaluate the action of $T^{(k)}(u)$  on symmetric functions. 
 Instead we use  (\ref{eq5_10}):
\begin{align*}
T^{(k)}(u)(f) &= \frac{1}{k+1}\partial ^{k+1}_u\st \,T(u,v) \st|_{v=u} (f)\\
&= \frac{1}{k+1}
\partial ^{k+1}_u \left(H(u) E(-v) E^{\perp}(-u) H^{\perp}(v)\right)|_{v=u} (f).
\end{align*}
Then the  Leibniz rule implies (\ref{eq6_2}).

\item 
By (\ref{eq5_6}), Corollary \ref{cor4_1} and the fact that  $H^\perp (v)(z^{r})=E^\perp (-u)(z^{r})=z^{r}$, one gets 
\begin{align}
H^\perp (v) S(u_l,\dots, u_1) ^{(m)} %\\
={\mathcal H_l(v^{-1})} S(u_l,\dots, u_1) ^{(m)},\label{eq6_5}
\end{align}
and
\begin{align*}
&H(u)E^{\perp}(-u)S(u_l,\dots, u_1) ^{(m)}\\
&=z^{m+l} u_l^m\dots u_1^m \prod_{1\le i<j\le l}\left( u_j- {u_i}\right)  H(u) \prod_{i=1}^{l} \left(1-\frac{u_i}{u}\right)
\;\prod_{i=1}^{l} H(u_i)\\
&={z^{-1} u^{-l-m}} S(u,u_l,\dots, u_1)^{(m)}.
\end{align*}
Substitution  of 
\begin{align*}
\frac{u^{m}}{v^{m}}&E(-v)H(u) E^{\perp}(-u)  H^{\perp}(v) S(u_l,\dots, u_1) ^{(m)} \\
&= \frac{\mathcal H_l(v^{-1}) E(-v)} { zu^l v^m}S(u, u_l,\dots, u_1)^{(m)}
\end{align*}
in (\ref{eq6_1})  proves (\ref{eq6_3}).
\item 
Differentiation of  the first equality in  Corollary \ref{cor4_1}  implies
  \[
 \partial_u^s(E^\perp(-u)) S(u_l,\dots, u_1)=  \partial_u^s \mathcal E_l(-u^{-1})S(u_l,\dots, u_1).
 \]
Using that $E(-u)=H(u)^{-1}$, $ {\mathcal E_l}(-u^{-1})= {\mathcal H_l}(u^{-1})$ and   (\ref{eq6_5}),
 we get  from (\ref{eq6_3}) 
    \begin{align*}
     T^{(k)}(u)& \,\left(S(u_l,.., u_1)\right)\\&=  \frac{1}{k+1}\sum_{s+r=k+1}{s\choose{ k+1}}  E(-u)  \partial^r_uH(u)  
     {{\mathcal H_l}(u^{-1})}  {\partial_u^s {\mathcal E_l}(-u^{-1}) }S(u_l,.., u_1)\\
   &= \frac{\,\partial_u^{k+1}({\mathcal E_l}(-u^{-1}) H(u)\,)}{(k+1) {\mathcal E_l}(-u^{-1}) H(u)} S(u_l,.., u_1).
    \end{align*}
\end{enumerate}
\end{proof}
\begin{remark}\label{rem6_1}%{rem6_3_1} 
Formula   (\ref{eq6_4}) carries remarkable properties, including  a  ``duality symmetry". 

 First, observe that  the action of $\hat a_\infty$ by (\ref{eq6_3}) and the action  of $\W_{1+\infty}$ by (\ref{eq6_4}) are expressed through  {\it multiplication operators} by symmetric functions 
 $h_k(x_1,x_2,\dots)$ and symmetric polynomials   $e_m(u_1,\dots u_l)$, and {\it contain no differentiation operators} (compare, for example,
  with (\ref{eq6_2}) that involves application of differential  operators $h_k^\perp (x_1,x_2,\dots)$ and $e_k^\perp (x_1,x_2,\dots)$).
  
 Second, (\ref{eq6_4}) interprets  formal distribution  $S(u_1,.., u_l)$ as an algebraic expression 
\[ 
S(u_1,.., u_l) \in \bC [[x_1,x_2,\dots]] \otimes \bC [[u_1, \dots u_l]],
\]
in other words,  composed  of products of  symmetric functions
 $s_\alpha= s_\lambda(x_1,x_2,\dots )\in \Lambda\subset \bC[[x_1,x_2,\dots]]$  and  monomials
   $u_1^{\alpha_1}\dots u_l^{\alpha_l}\in \bC[u_1, \dots u_l ]$, where  both components are treated on equal footing.
   In this interpretation the action of $T^{(k)}(u)$ reveals certain  duality: it is expressed through  
   action on the ``$\bC[[x_1,x_2,\dots]]$ part"   by {\it  multiplication} operators combined in   $ H(u; x_1,\dots, x_l)$,   and on the  ``$\bC[u_1, \dots u_l]$  part''   by action of  
    {\it multiplication} operators combined in 
    $\mathcal E_l(u; u_1,\dots, u_l)$. Note also that    operators  ``$H(u)$" and ``$\mathcal E_l(u)$"
  acting in different components can be considered of   ``dual nature"   due to relations of the kind (\ref{HEP1}).
The observed here combinatorial duality  are naturally connected to Remark \ref{rem5_2}.
\end{remark}

\subsection{Actions of Heisenberg and Virasoro algebra}\label{Sec6_2}
\begin{example} \label{ex6_1}
Formula (\ref{eq6_4})  interprets  the action of  Heisenberg algebra action on the basis  Schur symmetric functions  as a result of   multiplication of the following formal distributions:
\begin{align}\label{eq6_5}%{eqTs5}
  T^{(0)}(u) \,\left(S(u_1,.., u_l)\right)=&
 \left( \sum_{k\ge 1}{p_k} u^{k-1} +\sum_{i=1}^{l} i_{u_i/u}\left(\frac{u_i}{u(u-u_i)}\right)\right)\,S(u_1,.., u_l).
\end{align}
Indeed, from  (\ref{eq6_4}) and (\ref{HEP2})  we get
\begin{align*}
  T^{(0)}(u) \,\left(S(u_1,.., u_l)\right)&=\frac{\mathcal E_l(-u^{-1})  \partial_u H(u) + \partial_u \mathcal E_l (-u^{-1}) H(u)}
 {\mathcal E_l(-u^{-1})H(u)}\left(S(u_1,.., u_l)\right)\\
  &=\left(\frac{ \partial_u H(u)}{H(u)} +\frac{ \partial_u \mathcal E _l(-u^{-1})}
  {\mathcal E_l(-u^{-1})}\right) \left(S(u_1,.., u_l)\right)\\
  &= \left(P(u) + \frac{1}{u^2}{ \mathcal P_l}(u^{-1}) \right) S(u_1,.., u_l)\\
  &=\left( \sum_{k\ge 1}{p_k} u^{k-1} +  \sum_{k\ge 1}\sum_{ i=1}^{ l}\frac{u_i^k}{u^{k+1}}  \right)S(u_1,.., u_l).
\end{align*}
\end{example}

\begin{example} \label{ex6_2}
From  (\ref{eq6_4}),
(\ref{HEP3}),  and (\ref{HEP4}), we obtain 
the action of Virasoro algebra  on  Schur symmetric functions   in a form of  multiplication of the following formal distributions:
\begin{align*}
 & T^{(1)}(u) \,\left(S(u_1,.., u_l)\right)=\left( \frac{\partial^2_u H(u) }{2H(u)} + \frac{\partial_u\mathcal E_l(-u^{-1}) \,\partial_u H(u) }
 	 {\mathcal E_l(-u^{-1})  H(u) }+ \frac{ \partial^2_u \mathcal E_l(-u^{-1})}{2\mathcal E_l(-u^{-1}) }\right)S(u_1,.., u_l)\\
&=\left(\frac{\partial_u P(u)}{2}   +\frac{P(u)^2 }{2}  +   \, \frac{{\mathcal P_l }(u^{-1} )P(u)}{u^2}+  \frac{\partial_u \mathcal P_l(u^{-1}) }{2u^4} 
	+\frac{\mathcal P_l(u^{-1})^2 }{2u^4}  
	\right)
S(u_1,.., u_l)\\
&=\Bigl(\frac{1}{2} \partial_u P(u)  +\frac{1}{2} P(u)^2   
	+P(u)\sum_{i=1}^{l} i_{u_i/u}\left(\frac{u_i}{u(u-u_i)}\right)
	-\sum_{i=1}^{l} i_{u_i/u}\left(\frac{u_i}{2u^4(u-u_i)^2}\right) 
\\
&+\sum_{i,j=1}^{l} i_{u_iu_j/u} \left(\frac{u_i u_j}{2u^2(u-u_i)(u-u_j)}\right)
\Bigr)
S(u_1,.., u_l).
  \end{align*}
\end{example}
This is already complete information about action of Virasoro, Heisenberg algebras, encoded in relations on  generating functions. If for a certain application  one has to know 
the explicit action of the a  particular coefficient $T^{(0)}_k$  or $T^{(1)}_k$  on a particular Schur symmetric function, it  can be recovered in both Examples  \ref{ex6_1}  and \ref{ex6_2} by application of 
the Murnaghan–Nakayama rule, that states 
$
 p_{r}\cdot s_{\lambda }=\sum _{\mu }(-1)^{ht(\mu /\lambda )+1}s_{\mu },
$ 
where the sum is over all partitions $\mu$ such that  $\mu/\lambda$ is a rim-hook of size $r$ and $ht(\mu /\lambda )$ is the number of rows in the diagram $\mu/\lambda$,
(\cite{Stan}, 7.17).

\begin{example}
 Introduce  a formal distribution 
 \[
 J(u)=P(u)+P^\perp(u)=
 \sum_{k\ge 1}{ p_k}{ u^{k-1}}
 +
 \sum_{k\ge  1} k{\partial/}{\partial_{p_{k}}} {u^{-k-1}}.
\]  

 Formula  (\ref{eq6_7})  is 
 often used as    {\it the definition} of   the formal distribution of generators  $T^{(k)}(u)$. We deduce it here as {\it a corollary} of (\ref{eq6_2}).
\begin{proposition}
   	 \begin{align}\label{eq6_7}
   		 T^{(k)}(u)= \frac{1}{k+1}\st (\partial_u+J(u))^{k}(J(u))\st.
   	 \end{align}
 \end{proposition}
   \begin{remark}
   Formula  (\ref{eq6_7})  should be understood as  that the operator $ (\partial_u+J(u))^{k}$ is applied to the formal distribution  of multiplication operators $J(u)$, and then the resulting  formal distribution of operators  acting on the space of symmetric functions  is  put in the bosonic normal order.
    For example, 
    \[
    \st (\partial_u+J(u))(J(u))\st = \partial_u P(u)+ (P(u))^2+ 2P(u)P^\perp(u)+\partial_u P^\perp(u)+(P^\perp(u))^2.
    \]    
 \end{remark}   
      \begin{proof}
    We will prove (\ref{eq6_7}) by induction. 
      Base of induction:
     \[
      T^{(0)}(u)= P(u)+P^\perp(u)= J(u).
     \]
      Assume that  (\ref{eq6_7})   holds for some  $k$.  Then, by induction assumption, 
      \[
      \st (\partial_u+J)^{k+1}(J)\st =   \st (\partial_u+J)  (  \st (\partial_u+J)^{k}J  \st)\st=(k+1)\st (\partial_u+J) \left( T^{(k)}(u) \right)\st.
      \]

            By (\ref{eq6_2}), (\ref{HEP3}) and  (\ref{HEP4}),
        \begin{align*}
      T^{(k)}(u) =  \frac{1}{k+1}\sum_{
      \substack{ r+s=k+1\\r, s\ge 0}
      }{{k+1}\choose{s}} \,(\partial_u +P(u))^{r-1}(P(u) ) \,  (\partial_u +P^\perp(u))^{s-1}(P^\perp (u) ),
    \end{align*}
    where  for the uniform formula with some  minor abuse of notations we  assume $ (\partial_u +P^\perp(u))^{-1}(P^\perp (u) )=(\partial_u +P(u))^{-1}(P(u) ) =1$.
By Leibniz rule, 
     \begin{align*}      
    & \st (\partial_u+P+P^\perp)\, \left((\partial_u +P)^{r-1}(P ) \,  (\partial_u +P^\perp)^{s-1}(P^\perp )\right)\st\\
    & =      \st  \left((\partial_u +P)^{r}(P ) \,  (\partial_u +P^\perp)^{s-1}(P^\perp )\right)\st +
       \st  \left((\partial_u +P^{r-1}(P ) \,  (\partial_u +P^\perp)^{s}(P^\perp )\right)\st.
     \end{align*}  
          Hence, the coefficient of  the term    $ (\partial_u +P)^{r-1}(P ) \,  (\partial_u +P^\perp)^{s-1}(P^\perp )$
          in \\ $(k+1)\st (\partial_u+J)\, \left( T^{(k)}(u) \right)\st$  is ${{k+1}\choose {s}}+ {{k+1}\choose {s-1}}= {{k+2}\choose {s}}$, 
          and 
              \begin{align*}
               \st&(\partial_u+J)^{k+1}(J)\st =
                 (k+1)\st (\partial_u+J) \left( T^{(k)}(u) \right)\st \\
                 &=
                \sum_{\substack{r+s=k+2,\\ r, s\ge 0}}{{k+2}\choose{s}} \,(\partial_u +P(u))^{r-1}(P(u) ) \,  (\partial_u +P^\perp(u))^{s-1}(P^\perp (u) )
                = (k+2)  T^{(k+1)}(u).
              \end{align*}      
  \end{proof}
  \end{example}

  \subsection{Some connections to existing results}
As it is mentioned in the Introduction,   the explicit formulas for  the action of  operators from  $\W_{1+\infty}$  were studied by several authors. In particular, in  \cite{LY3}   explicit formulas for the action of coefficients  $T^{(k)}_r$ on  Schur symmetric functions $S_\mu$ were proved.   The authors of    \cite{LY3}  also  provided a list of references to the previous works that computed particular cases  of actions of such operators.  From those  we would like to mention that combinatorial  description   \cite{Baker} of the action of the Heisenberg algebra on symmetric functions is  the most close to  our approach, c.f. Example \ref{ex6_1}. Also we would like to outline the framework of formal distributions behind   the explicit formulas obtained in  \cite{LY3} connecting it with the set up of this note.  Namely,   by (\ref{t1}),
$
\partial_u^k T(u,v)= -\psi^{-}(v) \partial_u^k \psi^+(u) 
+\partial_u^k i_{u/v} \left(\frac{1}{v-u} \right),
$
and   a variation of calculation in the proof of (\ref{eq6_4}), followed by the expansion of  formal distributions, gives 
\begin{align*}
&\partial_u^k T(u,v) S(u_l,\dots, u_1)\\
=&
-\psi^{-}(v) \partial_u^k \psi^+(u)\psi^+(u_l)\dots \psi^+(u_1)(1)
+\partial_u^k i_{u/v} \left(\frac{1}{v-u} \right)\psi^+(u_l)\dots \psi^+(u_1)(1)\\
=&-\psi^{-}(v) \partial_u^k S(u,u_l,\dots, u_1)
+\partial_u^k i_{u/v} \left(\frac{1}{v-u} \right) S(u_l,\dots, u_1)\\
=& 
-\sum_{s,a\in \bZ,\, \mu\in \bZ^l}  (a)_{k} \psi^- _{s-1/2}(S_{(\mu, a)} )v^{-s}u_1^{\mu_1}\dots u_l^{\mu_l}u^{a-k}
 \\
&+
\sum_{a=0}^\infty \sum_{ \mu\in \bZ^l} ({a)}_k S_{\mu}u_1^{\mu_1}\dots u_l^{\mu_l} {u^{a-k}}{v^{-a-1}}\\
=&\sum_{s,a\in \bZ,\, \mu\in \bZ^l} 
(a)_{k}
\left(
- \psi^- _{s-1/2}(S_{\mu, a} )
+\delta_{a\ge0}\delta_{a, s+1} S_{\mu}
\right)u_1^{\mu_1}\dots u_l^{\mu_l}u^{a-k}v^{-s}.
\end{align*}
Evaluating at $v=u$ we obtain
\begin{align*}
 &T(u)^{(k)}S(u_l,\dots, u_1)\\
 &=
 \sum_{r\in \bZ,\,  \mu\in \bZ^l} \sum_{s\in \bZ} 
(r+s)_{k}
\left(
- \psi^- _{s-1/2} (S_{\mu, r+s })
+\delta_{r+s \ge 0}\delta_{r, 1} S_{\mu}
\right) 
u^{r-k }u_1^{\mu_1}\dots u_l^{\mu_l}
\\
&=
 \sum_{r\in \bZ,\,  \mu\in \bZ^l}\left( -\sum_{s\in \bZ} 
(r+s)_{k}\psi^- _{s-1/2} (S_{\mu, r+s })
+ \delta_{r, 1} \sum_{s\ge -1} (s+1)_{k} S_{\mu}
\right) 
u^{r-k }u_1^{\mu_1}\dots u_l^{\mu_l}.
\end{align*}
Substitution of the combinatorial expression for the value of the action $\psi^- _{s-1/2} (S_{\mu, r+s })$ of charged free fermions  on the Schur symmetric function  in this formula gives the value of the coefficient  of the monomial  $u^{r-k }u_1^{\mu_1}\dots u_l^{\mu_l}$ and describes the action of $T^{(k)}_r$ on $S_\mu$ by a combinatorial formula  equivalent to \cite{LY3},  
where all calculations  are carefully performed and the final result is provided.

\section{Lie algebra $\W^B_{1+\infty}$ as a subalgebra of $ \hat a_\infty$} \label{Sec7}
\subsection{Lie algebra  $ \hat o_\infty$} \label{Sec7_1}

The B-type analogue of $\W_{1+\infty}$  was studied in \cite { A2, LY3, VDL}, and we follow these papers for definitions. 
Define an  anti-involution $i$  on $ a_\infty$ by its action  on generators 
\[
i( E_{kl})= (-1)^{k+l}E_{-l,-k}, 
\]
or, equivalently, by the   transformation of the matrix of generators
\[
i(T(u,w))= - \frac{w}{u}T(-w,-u).
\]
Let 
\begin{align}\label{eq7_1}%{eq10_1}
T^B(u,w)= uT(u,-w)-wT(w,-u)=  u(T(u,-w) -i(T(u,-w)).
\end{align}
Then 
\begin{align*}
i(T^B(u,w))= -T^B(u,w) =T^B(w,u),\quad \text{and}\quad T^B(u,u)=0.
\end{align*}
Coefficients of expansion of $T^B(u,w)$  are given by 
\[
T^B(u,w)= \sum_{i,j\in \bZ} (-1)^j\hat F_{ij}u^i w^{-j},\quad \text{where} \quad \hat F_{ij}= \hat E_{ij}-(-1)^{i+j} \hat E_{-j \,-i}.
\]

  \begin{proposition}\label{prop9_1} 
   $[ T^B(u,w), C]=0$, and 
  \begin{align*}
  [ T^B (u,w),  T^B (v,z)]&=
  2v\delta(v,-w) \left(T^B (u,z) + i_{z/u}\left(\frac{2u}{u+z}\right)C\right)\\
  &-2v\delta(v,-u) \left(T^B (w,z) + i_{z/w}\left(\frac{2w}{w+z}\right)C\right)\\
 & +\,2z\delta(z,-w) \left(T^B (v,u) +i_{u/v}\left(\frac{2v}{v+u}\right)C\right)\\
 &-2z\delta(z,-u) \left(T^B (v,w) +i_{w/v}\left(\frac{2v}{v+w}\right)C\right).
  \end{align*}
    \end{proposition}
  \begin{proof}
Statement follows from (\ref{eq7_1}) and  commutation relations  (\ref{eq2_4}) in  $\hat a_\infty$.
      \end{proof}
      
       \begin{remark}\label{rem9_1}
  The central part of this commutation relation can be expressed through $\gamma(u,w,v,z)$  defined by (\ref{eq2_5})
  as  
   \begin{align*}
  &2v\delta(v,-w)  i_{z/u}\left(\frac{2u}{u+z}\right)
  -2v\delta(v,-u)  i_{z/w}\left(\frac{2w}{w+z}\right)\\
 &+\,2z\delta(z,-w) i_{u/v}\left(\frac{2v}{v+u}\right)
 -2z\delta(z,-u)i_{w/v}\left(\frac{2v}{v+w}\right)\\
 &\quad = 4vu\gamma(u,-w,v-z)+4vw\gamma(-w,-u,v, z).
  \end{align*}
  \end{remark}
 
  Let $o_\infty=\{a\in a_\infty|  i(a)=-a\}$, and let  $\hat o_\infty= o_\infty\oplus \bC C$. Together with commutativity of the  central element $C$,   Proposition \ref{prop9_1}  expresses commutation relations between generators $\{\hat F_{ij}\}$  of $ \hat o_\infty$ in terms of relations of generating functions $T^B(u,w)$.
 
  \subsection{Inifinite-dimensional Lie algebra  $ \W^B_{1+\infty}$} \label{Sec7_2}
  Let $\W_{1+\infty}^B=\W_{1+\infty}\cap \hat o_\infty $.  As it is observed in \cite{VDL},  every element of  $\W_{1+\infty}^B$ is  a linear span of  expressions of the form 
$\sum_{j\in \bZ}f(-j)F_{k-j, j}$ with polynomial coefficients  $f(-j)$. Since 
\begin{align*}
i(T^{(k)}_r)= (-1)^r\sum_{j\in \bZ}(r+j-1)_kE_{-j,-r-j}=(-1)^r\sum_{j\in\bZ}(-j-1)_kE_{j+r,-j}, 
\end{align*}
this implies that $i(T^{(k)}_r) \in \W_{1+\infty}$, and   $T^{(k)}_r-i(T^{(k)}_r)\in \W_{1+\infty}^B$. Note also
\[
T^{(k)}(u)-i(T^{(k)}(u))= \sum_r \left(T^{(k)}_r-i(T^{(k)}_r)\right)u^{k-r}= \sum_r (r+j-1)_kF_{r+j,j}u^{k-r}.
\]
On the basis elements of  the algebra $\hat D$, the  defined above anti-involution has the form
\[
i(t^kD^l)= (-D)^l (-t)^k= (- t)^k (-D-k)^l.
\]
  
Following  \cite {A2, LY3}, define 
 \[
T^{B(k)}(u)=\partial^k_xT^B(x,y)|_{x=u, y=-u}.
\]
In the expansion 
$
T^{B(k)}(u)=\sum_{r\in \bZ} T^{B(k)}_r u^{r-k}
$
the coefficients have the form\\ $T^{B(k)}_r=\sum_{j\in \bZ}(r+j)_kF_{j+r,j}$,
and  using 
that
\[
(a+1)_k= (a)_k+k(a)_{k-1},\quad
(a-1)_k= (a)_k-k(a)_{k-1},
\]
we get relation between $T^{B(k)}(u)$ and $T^{(k)}(u)$:
\begin{align*}
T^{B(k)}(u)&=T^{(k)}(u)-i(T^{(k)}(u)) + k(T^{(k-1)}(u)-i(T^{(k-1)}(u))), \quad 
\\
T^{(k)}(u)&-i(T^{(k)}(u)) = T^{B(k)}(u)- kT^{B(k-1)}(u). \quad 
\end{align*}
In particular,  these relations show that  $T^{B(k)}_r\in \W_{1+\infty}^B$.
\begin{example} For 
$k=0$ one gets
$
T^{B(0)}(u)= u \,T^{(0)}(u)-  uT^{(0)}(-u),\,
$
and from  Proposition \ref{prop9_1},  
\[
[T^{B(0)}(u), T^{B(0)}(v)]=  
 4uv\left(\partial_v\delta(u,v)-\partial_v\delta(-u,v)\right).
\]

\end{example}

%%%%%%%%%%%%%%%%%%%%%%%%%%%%%%
\section{Presentation  of $\hat o_\infty$ and  $\W^B_{1+\infty}$ through neutral  fermions}\label{Sec8}
Similarly to  Section \ref{Sec3}, we will  use the action of a Clifford algebra to construct the action of generators of $\hat o_\infty$ and $\W^{B}_{1+\infty}$
on symmetric functions. 
\subsection{Clifford algebra of neutral  fermions}\label{Sec8_1}

\begin{definition}\label{def8_1}
The Clifford algebra of {\it neutral fermions} is  generated by  
 $\{\varphi_i\}_{i\in \bZ}$  with relations 
 \begin{align}\label{eq8_1}
  \varphi_m \varphi_n+\varphi_n \varphi_m= 2(-1)^m \delta_{m+n,0}\quad \text {for}\quad m, n\in \bZ.
  \end{align}
  \end{definition}
Collecting  generators as coefficients  of a formal distribution  
 $\varphi(u)=\sum_{j\in \bZ}\varphi_j u^{-j}$,  relations (\ref{eq8_1}) can be written as
\[
 \varphi(u)\varphi(v) + \varphi(v)\varphi(u) =2v\delta(v,-u).
 \]
With  the same Definition  \ref{def3_1} of normal ordered product,  by direct calculation 
we get the following  relations.
\begin{lemma}\label{lem8_1}
\begin{eqnarray*}
  :\varphi(u)\varphi(v):&= \varphi (u) \varphi(v) -2 i_{v/u} \left(\frac{u}{u+v} \right)=-\varphi(v)\varphi(u)+2i_{u/v}\left(\frac{v}{v+u}\right).
  \end{eqnarray*}
  \end{lemma}

\subsection{Presentation  through neutral fermions}\label{Sec8_2}

\begin{proposition} \label{prop8_1}
The algebra $\hat o_\infty$ can be realized through  neutral fermions  by the identification
 \begin{align*}
  T^{B}(u,v)\to :\varphi(u)\varphi(v): \quad \text{ and $C\to Id$.}
  \end{align*}
 Accordingly, 
  \begin{align*}
  T^{B(k)}(u,v)\to :\partial^k_u\varphi(u)\varphi(v): .
  \end{align*}
\end{proposition}

  \begin{proof}
  Using Lemma \ref{lem8_1},  we  express  commutators of  normal order products through commutators of regular products, rearrange the order of  terms.
   \begin{align*}
    &[:{\varphi}(u){\varphi}(w):,:{\varphi}(v){\varphi}(z):]
   =[{\varphi}(u){\varphi}(w),{\varphi}(v){\varphi}(z)]\\
    &={\varphi}(u){\varphi}(w){\varphi}(v){\varphi}(z)-{\varphi}(v){\varphi}(z){\varphi}(u){\varphi}(w)\\
    &=2v\delta(v,-w){\varphi}(u){\varphi}(z)-2v\delta(v,-u){\varphi}(w){\varphi}(z)
    +2z\delta(z,-w){\varphi}(v){\varphi}(u)-2z\delta(z,-u){\varphi}(v){\varphi}(w).
\end{align*}
Again with the help of Lemma  \ref{lem8_1},  we  express  regular products  through normal order products to  conclude that $:\varphi(u)\varphi(v):$ satisfies exactly 
the same commutation relation as $T^{B}(u,v)$ in Proposition \ref{prop9_1}.
  \end{proof}
  
  \begin{proposition} \label{prop8_2} 
 \begin{align*}
[T^{B}(u,v),\varphi(z)]&= 2z\delta(z,-u) \varphi(v)-2z\delta(z,-v) \varphi(u),
\\
[T^{B(k)}(u),\varphi(z)]&=-2z\partial_u^k\delta(u,-z) \varphi(-u)-2z\delta(z,u) \partial_u^k\varphi(u).
  \end{align*}
\end{proposition}
    \begin{proof}
    By Proposition \ref{prop8_1} and Lemma \ref{lem8_1},
    \begin{align*}
      [T^{B}(u,v),\varphi(z)]&=[
      \varphi(u)\varphi(v),\varphi(z)]=2z\delta(z,-u) \varphi(v)-2z\delta(z,-v) \varphi(u).
    \end{align*}
    Using (\ref{eq2_21}) and differentiating both sides of the line above,  we get the second statement. 
  \end{proof}

\section{Review of Properties of Schur $Q$-functions}\label{Sec9}
 The analogue of  Theorem \ref{thm6_1} that would describe the action of generators of $\hat o_\infty$ and $\W^B_{1+\infty}$ on  particular families of symmetric functions can be obtained by direct calculations. However,  below,   whenever it is possible,   we   aim   {\it to reduce our argument} to the work already done  for type $A$, instead of repeating it     in  the complete analogy to  type $A$.
With this goal in  mind, in this section  we review how  Schur  $Q$-functions are connected to Schur symmetric functions. 
\subsection{ Symmetric Schur  $Q$-functions }\label{9_1} 
Let $q_k= q_k(x_1,x_2,\dots)$ be symmetric functions in variables $x_1,x_2, \dots$, defined  as coefficients of the expansion of formal series  $Q(u)\in  \Lambda[[u]]$, where
 \begin{align}\label{eq9_1}%{shurq}
Q(u) =\sum_{k\in \bZ} q_k u^k= E(u) H(u)=\prod_{i\in \bN} \frac{1+x_iu}{1-x_iu}.
\end{align}
Then $q_k=\sum_{i=0}^k e_ih_{k-i}$ for $k>0$, $q_0=1$, and $q_k=0$ for $k<0$.

\begin{proposition}\label{prop_rel2}  \cite[III.8]{Md} We have in  $\mathcal D[[u]]$ (respectively,  in  $\mathcal D[[ u^{-1}]]$\,) 
  \begin{align*}
Q(u) =S_{odd}(u)^2, \quad\text{where}\quad   
S_{odd}(u)=exp\left(\sum_{n\in \bN_{odd}}\frac{p_{n}}{n}{u^{n}}\right),
\end{align*}
\begin{align*}
S_{odd}^{\perp}(u)=exp\left(\sum_{n\in \bN_{odd}} \frac{\partial}{\partial p_n}\frac{1}{u^{n}}\right),
\end{align*}
where $\bN_{odd}=\{1,3,5,\dots\}$.
\end{proposition}

%%%%%%%%%%%%%%%%%%
Consider formal expansion of the rational function
 \begin{align*}
% f(u,v)=
 i_{v/u}\left(\frac{u-v}{u+v}\right)=\left(1-\frac{v}{u}\right)\sum_{k\in \bZ_+}(-1)^k\frac{v^k}{u^k} =1+2\sum_{k\in \bN}(-1)^k\frac{v^k}{u^k}\in \bC[[v/u]].
 \end{align*}
 Note that 
 \begin{align}\label{eq9_1}%{fuv}
  i_{v/u}\left(\frac{u-v}{u+v}\right)+  i_{u/v}\left(\frac{v-u}{v+u}\right)
=(u-v)\delta(u,-v) =2\sum_{k\in \bZ}(-1)^k\frac{u^k}{v^k}.
 \end{align} 
Then  {\it Schur Q-function }
$Q_\lambda =Q_\lambda (x_1, x_2, \dots)$ can be defined  as  coefficients of ${u^{\lambda_1}\dots u^{\lambda_l}}$
of the formal distribution
\begin{align}\label{eq9_3}%{genQ}
Q(u_l,\dots, u_1)=\sum_{\lambda_1,\dots, \lambda_l\in \bZ} {Q_\lambda }{u_1^{\lambda_1}\dots u_l^{\lambda_l}}=\prod_{1\le i<j\le l}i_{u_j/u_i}\left(\frac{u_i-u_j}{u_i+u_j} \right)\prod_{i=1}^{l} Q(u_i),
\end{align}
see \cite[III, (8.8)]{Md}. %CHECKED AUG 30 2025
 The values of these coefficients as symmetric polynomials 
 are  given by the formula
\begin{align}\label{eq9_4}%{defQ}
Q_\lambda (x_1,\dots, x_N)=
\frac{2^l}{(N-l)!}
\sum_{\sigma\in S_N} \prod_{i=1}^{l}x_{\sigma(i)}^{\lambda_i} \prod_{i<j}\frac{x_{\sigma(i)}+x_{\sigma(j)}}{x_{\sigma(i)}- x_{\sigma(j)}},
\end{align}%CHECKED AUG 30 2025
for strict partitions $\lambda$, and are zero  otherwise
  \cite[III, (8.7)]{Md}. %CHECKED AUG 30 2025
Schur $Q$-polynomials have a stabilization  property \cite[III, (2.5)]{Md},  hence  one can omit the number $N$ of variables $x_i$'s, as long as it is  not less than the length of the partition $\lambda$, considering  $Q_\lambda$ as a function of  infinitely many variables $(x_1, x_2 \dots )$.

\subsection{Action of neutral fermions  on   the bosonic space $\B_{odd}$  }\label{9_2} %%%%%%%%% Section Q

Let  $\B_{odd}$ be the  subspace of the 
ring   of symmetric functions generated by odd power sums:
  \[\B_{odd}= \bC[  p_1, p_3, p_5,\dots]\subset \Lambda=\bC[  p_1, p_2, p_3,p_4\dots].
   \]
   Then   $q_k\in \B_{odd}$,  and  that $\B_{odd}= \bC[q_1,q_3,\dots ]$,  \cite[ III.8 (8.3)]{Md}.

It is not difficult  to prove  (see e.g. \cite{NRQ})  that the space  $\B_{odd}$  is invariant with respect  to action of  $e^\perp_k$,  and $ h^\perp_k$, and 
that the  restrictions to  ${\B_{odd}}$ of action of coefficients of the following  formal distributions  coincide:
\[
 E^{\perp}(u)|_{\B_{odd}}= H^{\perp}(u)|_{\B_{odd}}=S_{odd}^{\perp}(u).
 \]
Denote as $ E_{odd}^{\perp}(u)=  E^{\perp}(u)|_{\B_{odd}}$ and  $ H_{odd}^{\perp}(u)=  H^{\perp}(u)|_{\B_{odd}}$.

Define a  formal distribution  $\varphi(u)$ of operators acting on $\B_{odd}$:
 \begin{align}\label{eq9_5}%{defphi1}
&\varphi(u)=  E(u)H(u) E^{\perp}_{odd}(-u)=Q(u) S_{odd}^{\perp}(-u).
\end{align}

Let  $\{\varphi_i\}_{i\in \bZ}$ be
coefficients  of the expansion $
\varphi(u)=\sum_{j\in \bZ}\varphi_{-j} u^{j}$.
One  can  prove  (see e.g. \cite{NRQ})  that   $\varphi(u)$ is a quantum field that  acts exactly as   the Clifford algebra of {\it neutral fermions} on the space $\B_{odd}$, satisfying  relations
\[
 \varphi(u)\varphi(v) + \varphi(v)\varphi(u) =2v\delta(v,-u),
 \]
 or, in terms of coefficients, 
 \begin{align}\label{eq9_6}%{neut1}
\varphi_m \varphi_n+\varphi_n \varphi_m= 2(-1)^m \delta_{m+n,0}\quad \text {for}\quad m, n\in \bZ.
\end{align}
This quantum field can be also expressed in the bosonic form     
\begin{align*}%\label{phiqr}
\varphi(u)=  Q(u) S_{odd}(-u)^{\perp}= exp\left(\sum_{n\in \bN_{odd}}\frac{2p_{n}}{n}{u^{n}}\right)exp\left(-\sum_{n\in \bN_{odd}} \frac{\partial}{\partial p_n}\frac{1}{u^{n}}\right).
\end{align*}
Note that  (\ref{eq9_5}) implies that the action of neutral fermions on $\B_{odd}$ can be expressed through the restrictions  of the action of charged free fermions on $\B_{odd}$:
 \begin{align*}%\label{defphi2}
&\varphi(u)= H(u)\psi^+(u)|_{\B_{odd} }=  E(u)\psi^-(-u)|_{\B_{odd} }.
\end{align*}

By  (\ref{eq9_5}) and  Lemma \ref{lemmaHE} we get % Vrode verno  Aug 11 2025
variations of  vertex operator presentations of Schur $Q$-functions, first considered in 
  \cite{Jing3}:
\begin{align}%label{eq9_6}%{defQBKP}
\varphi(u_l)\dots \varphi(u_1) (1)&= Q(u_l,\dots, u_1), \label{eq9_7}\\
\varphi_{-\lambda_l}\dots \varphi_{-\lambda_1} (1)&= Q_\lambda.\notag
\end{align}

\subsection{Bosonic normal ordered products for neutral fermions}\label{9_3}

Similarly to Section \ref{Sec5_3}, define the bosonic normal ordered product of neutral fermions as 
\begin{align}\label{eq9_8}
 \st \varphi(u)\varphi(v)\st=Q(u)Q(v) E_{odd}^{\perp}(-u)E_{odd}^{\perp}(-v).
\end{align}
 
Relations between different products  of neutral free fermions are collected in the following lemma that  follows  from  the definition of $\st \varphi(u)\varphi(v)\st$ and  Lemma \ref{lemmaHE}.
\begin{lemma}\label{lem9_1}
  \begin{align}\label{eq9_9}
    :\varphi(u)\varphi(v):&= \varphi (u) \varphi(v) -2 i_{v/u} \left(\frac{u}{u+v} \right)=
       i_{v/u} \left( \frac{u-v}{u+v}\right) \st \varphi(u)\varphi(v)\st  -     2i_{v/u} \left( \frac{u}{u+v}\right).
\end{align}
\end{lemma}
 %%%%%%%%% Section Q
We introduce notation 
\[
\st T^B(u,v)\st = \st \varphi(u)\varphi(v)\st.
\]
Note that this formal distribution  can be expressed through restriction of formal distribution $\st T(u,v)\st$ to $\B_{odd}$:
\[
\st T^B(u,v)\st =  E(u)H(v)  \, \st T(u,-v)\st |_{\B_{odd}},
\]
and that 
$
\st T^B(u,-u)\st = 1.
$

Define a  differential operator
\[
D^{(k)}_u= \partial_u^k+ \frac{2u}{k+1}\partial_u^{k+1}.
\]
\begin{lemma} \label{lem9_2}%{lem12_2}
Operator  $D^{(k)}_u$ satisfies the analogue of Leibniz rule: 
\[
D^{(k)}_u(f g)= \sum_{r=0}^k{k\choose r}\partial_u^r (f)\,\, D^{(k-r)}_u( g) + \frac{2u}{k+1}\partial^{k+1}_u( f)\,  g.
\]
\end{lemma}

\begin{proposition} \label{prop9_2}%prop_122}
 We can express  the generating functions of $\hat o_\infty$ and  of $\W^B_{1+\infty}$ as
\[
 T^B(u,v)=i_{v/u} \left( \frac{1}{u+v}\right) 
  \left(  (u-v)\st\, T^B(u,v)\st -2u\right),
\] % Caleb confirmed April 8 25
and
\begin{align}\label{eq9_11}%eq_134}
 T^{B\, (k)} (u)&=D^{(k)}_u \st\, T^B(u,v))\st|_{v=-u}  - 2\delta_{k,0}.
\end{align}% Caleb confirmed April 8 25
\end{proposition}
\begin{proof}
The first statement immediately follows from Proposition \ref{prop8_1} and  Lemma \ref{lem9_1}.
For the second statement,
observe that 
\[
(u+v) T^B(u,v)=
 (u-v)\st\, T^B(u,v)\st -2u.
\] Taking the derivative $\partial_u^{k+1}$ of both sides of this equality   and evaluating $v=-u$ gives
\begin{align*}
 (k+1)T^{B(k)}(u)=
(k+1)\partial_u^{k}\st\,T^B(u,v)\st\, |_{v=-u}+ 2u \partial_u^{k+1}\st\,T^B(u,v)\st|_{v=-u}\,
 -2\delta_{k,0},
\end{align*}
which proves (\ref{eq9_11}).
\end{proof}

\section{Action of $\hat o_\infty$ and $\W^B_{1+\infty}$ on  symmetric functions}\label{Sec10}
\subsection{The second main result}
Introduce notation
\begin{align*}
\mathcal Q_l(-u^{-1})=\prod_{s=1}^{l} i_{u_i/u}\left(\frac{1+u_i/u}{1-u_i/u}\right)= \mathcal E_l(-u^{-1}) \mathcal H_l(-u^{-1}).
\end{align*}

\begin{theorem}\label{thm10_1}
 	 \begin{enumerate}
		  \item  The action of $\hat o_\infty $  on  any  $f\in \B_{odd}$  is described by 
			\begin{align}\label{eq10_1}%{Tq1}
				T^B(u,v)(f)=\left(i_{v/u}\left(\frac{u-v}{u+v}\right) Q(u) Q(v) E^\perp_{odd}(-u)E^\perp_{odd}(-v) - i_{v/u}\left(\frac{2u}{u+v}\right)\right) (f).
			\end{align}
		\item  The action of $\hat o_\infty $ on the basis of Schur symmetric $Q$-functions  in terms of generating functions is given by
			\begin{align}\label{eq10_2}%{Tq2}
				T^B(u,v)&\left( Q(u_l,.., u_1)\right)
			=\quad  Q(u,v,u_l, \dots, u_1) -  i_{v/u}\left(\frac{2u}{u+v}\right) Q(u_l, \dots, u_1).
			\end{align}
		\item  The action of  $\W^B_{1+\infty}$ on  any $f\in \B_{odd}$  is described by  
 				 \begin{align}\label{eq10_3}%{Tq3}
   					 T^{B\,(k)}(u) \,(f)= &  
					\left(\sum_{r=0}^k{k\choose r}\partial_u^r( Q(u))Q(-u)\,\, E_{odd}^{\perp}(u) D^{(k-r)}_u (E_{odd}^{\perp}(-u))\right)(f)\\
    					&+\left(\frac{2u}{k+1}\partial^{k+1}_u Q(u) Q(-u) - 2\delta_{k,0}\right)(f).\notag
 				 \end{align}
		\item The action of  $\W^B_{1+\infty}$ on the basis of Schur symmetric $Q$-functions is  defined by
				 \begin{align} \label{eq10_4}%{Tq4}
  					 T^{B\,(k)}(u) \,\left(Q(u_l,.., u_1)\right)
   					= D^{(k)}_u (Q(u) \mathcal Q_l(u))\,Q(-u, u_l,.., u_1)- 2\delta_{k,0}Q(u_l,.., u_1).
 				 \end{align}
		\end{enumerate}
\end{theorem}

%%%%%%%%%%%%%%%%Proof
\begin{proof}
\begin{enumerate}
\item  %%%%%
Follows immediately from (\ref{eq9_8}) and  Proposition \ref {prop9_2}.
\item %%%%%
By Lemma \ref{lem9_1} and  (\ref{eq9_7}),
 \begin{align*}
 T^B(u,v)&\left( Q(u_l,.., u_1)\right)= \left(\varphi(u)\varphi(v) - i_{v/u}\left(\frac{2u}{u+v}\right)\right)\left( Q(u_l,.., u_1)\right)\\
 = &Q(u,v,u_l,.., u_1) -i_{v/u}\left(\frac{2u}{u+v}\right) Q(u_l,.., u_1).
  \end{align*}
\item %%%%%
Using Lemma \ref{lem9_2},  and that $E_{odd}^{\perp}(-u)E_{odd}^{\perp}(u)= 1$, for $k>0$ we get
\begin{align*}
&T^{B\,(k)}(u)=D^{(k)}_u \st\, T^B(u,v))\st|_{v=-u}  - 2\delta_{k,0}\\
&=D^{(k)}_u \left(Q(u) Q(v)E_{odd}^{\perp}(-v)E_{odd}^{\perp}(-u)\right)\st|_{v=-u}  - 2\delta_{k,0} \\
&= \sum_{r=0}^k{k\choose r}\partial_u^r (Q(u))Q(-u)\,\,E_{odd}^{\perp}(u) D^{(k-r)}_u (E_{odd}^{\perp}(-u))+ \frac{2u}{k+1}\partial^{k+1}_u Q(u) Q(-u)  - 2\delta_{k,0}.
\end{align*}
\item %%%%
Note that
\[
E^\perp_{odd}(-u) Q(u_l,\dots, u_1)=\mathcal Q_l(-u^{-1})Q(u_l,\dots, u_1),
\]
and  therefore
\[
D^{(k-r)}_u(E^\perp_{odd}(-u) Q(u_l,\dots, u_1))=D^{(k-r)}_u(\mathcal Q_l(-u^{-1})) Q(u_l,\dots, u_1).
\]
From (\ref {eq10_3}), (\ref{eq9_5}), and Lemma \ref {lem9_2},  for $k>0$,
\begin{align*}
&T^{B\,(k)}(u)Q(u_l,\dots, u_1)\\
&= \sum_{r=0}^k{k\choose r}\partial_u^r (Q(u))Q(-u)\,\,E_{odd}^{\perp}(u) D^{(k-r)}_u(\mathcal Q_l(-u^{-1})) Q(u_l,\dots, u_1)\\
&+ \frac{2u}{k+1}\partial^{k+1}_u Q(u) Q(-u)Q(u_l,\dots, u_1).\\
&= \sum_{r=0}^k{k\choose r}\partial_u^r (Q(u)) D^{(k-r)}_u(\mathcal Q_l(u)) \, Q(-u, u_l,\dots, u_1)\\
&+ \frac{2u}{k+1}\partial^{k+1}_u Q(u) Q(-u)Q(u_l,\dots, u_1).\\
&=  D^{(k)}_u(Q(u) \mathcal Q_l(u)) Q(-u, u_l,\dots, u_1)- \frac{2u}{k+1}\partial^{k+1}_u (Q(u))\, \mathcal Q_l(u)Q(-u, u_l,\dots, u_1)\\
& + \frac{2u}{k+1}\partial^{k+1}_u Q(u) Q(-u)Q(u_l,\dots, u_1).
\end{align*}
Observe that,  by (\ref{eq9_3}),
$Q(-u, u_l,\dots, u_1)= Q(-u) \mathcal Q_l(-u) Q(u_l, \dots, u_1),
$
and that $\mathcal Q_l(u)\mathcal Q_l(-u)=1$.
This implies that 
\begin{align*}
&- \frac{2u}{k+1}\partial^{k+1}_u Q(u) \mathcal Q_l(u)Q(-u, u_l,\dots, u_1)
 + \frac{2u}{k+1}\partial^{k+1}_u Q(u) Q(-u)Q(u_l,\dots, u_1)=0.\\
\end{align*}
For $k=0$ one more term of the form $2\delta_{k,0}Q(u_l,.., u_1)$ is subtracted, which completes the proof of proposition. 
\end{enumerate}

\end{proof}
\begin{remark}\label{rem10_1}
In analogy to Remark \ref{rem6_1},  we observe that the action (\ref{eq10_3}) of $\hat o_\infty$ and the action  (\ref{eq10_4}) of $\W^B_{1+\infty}$ are again described by multiplication operators, with  a  duality property naturally connected to the Cauchy identity  for Schur $Q$-functions
\begin{align*}
Q(u_l,\dots, u_1)\prod_{1\le i<j\le l}i_{u_j/u_i}\left(\frac{u_i+u_j}{u_i-u_j} \right)=\prod_{i=1}^{l} Q(u_i)
=
\sum_{\lambda\in \bZ}  {Q_\lambda(x_1,x_2,\dots)} Q_\lambda({u_1,\dots, u_l}),
\end{align*}
for $|u_1|<|u_2|<\dots <|u_l|$, \cite[III.8 (8.13)]{Md}.
\end{remark}

\subsection*{Acknowledgments}The research of the
second author is supported by the Simons Foundation Travel Support for Mathematicians grant MP-TSM-00002544.

\end{document}